\documentclass[12pt,reqno]{amsart}
\usepackage{amsmath,amsfonts,amssymb,amsthm,cite,color}
\usepackage[english]{babel}
\usepackage[left=2.9cm,right=2.9cm,top=2.8cm,bottom=2.8cm]{geometry}
\usepackage[colorlinks=false,pdfborder={0 0 0},linktocpage,
pdfpagelabels,bookmarksnumbered,bookmarksopen]{hyperref}
\numberwithin{equation}{section}
\newtheorem{theo}{Theorem}[section]
\newtheorem{Lemma}{Lemma}[section]
\newtheorem{Proposition}{Proposition}[section]
\theoremstyle{remark}
\newtheorem{Remark}{Remark}[section]
\theoremstyle{definition}
\newtheorem{Definition}{Definition}[section]
\title[Solutions concentrating on spheres]{Solutions concentrating on spheres for Hamiltonian elliptic systems with supercritical growth}

\author[H. Zhang]{Hui Zhang}
\author[J. J. Zhang]{Jianjun Zhang}
\author[X. X. Zhong]{Xuexiu Zhong}

\address[H.\ Zhang]{\newline\indent Department of  Mathematics,
Jinling Institute of Technology,
\newline\indent
Nanjing 211169, PR China}
\email{\href{mailto:huihz0517@126.com}{huihz0517@126.com}}

\address[J. J. \ Zhang]{\newline\indent College of Mathematics and Statistics, Chongqing Jiaotong University,
\newline\indent
Chongqing 400074, PR China}
\email{\href{mailto:zhangjianjun09@tsinghua.org.cn}{zhangjianjun09@tsinghua.org.cn}}

\address[X.~X.~Zhong]{\newline\indent South China Research Center for Applied Mathematics and Interdisciplinary Studies,
\newline\indent
South China Normal University
\newline\indent
Guangzhou 510631, China}
\email{\href{mailto:zhongxuexiu1989@163.com}{zhongxuexiu1989@163.com}}

\subjclass[2000]{35J50; 35B25; 35B06.}
\keywords{Hamiltonian elliptic system; Symmetry; Concentration phenomena; Generalized Nehari manifold.}

\begin{document}
\begin{abstract}

We study positive solutions of the singularly perturbed Hamiltonian elliptic system
\begin{equation*}
\begin{cases}
-\varepsilon^2\Delta u+u=|v|^{q-2}v &\text{in }A,\\
-\varepsilon^2\Delta v+v=|u|^{p-2}u &\text{in }A,\\
u,v>0\text{ in }A,\qquad u=v=0 &\text{on }\partial A,
\end{cases}
\end{equation*}
where $A=\{x\in\mathbb R^{2m}:a<|x|<b\}$, $0<a<b$, $m\geq2$, and
$2<p,q<2(m+1)/(m-1)$. We prove the existence of a family of solutions whose two
components concentrate around the same $(m-1)$-dimensional sphere on the inner boundary
as $\varepsilon\to0$. The two spheres of maxima are at distance $o(\varepsilon)$
from each other, and both components decay exponentially away from these spheres.
The admissible range of $(p,q)$ includes critical and supercritical exponent pairs in
dimension $2m$. The proof uses a reduction to a weighted system in dimension
$m+1$. We establish foliated Schwarz symmetry of its ground states with a common
axis via reduction by inversion and polarization. A uniform estimate on generalized
Nehari fibers gives the upper energy bound. A half-space nonexistence result and
the corresponding lower bound then determine the concentration location and the
limiting profiles.

\end{abstract}
\maketitle

\section{Introduction and main results}
In this paper, we consider the singularly perturbed Hamiltonian elliptic system
\begin{equation}\label{1.1}
\begin{cases}
-\varepsilon^2\Delta u+u=|v|^{q-2}v &\text{in }A,\\
-\varepsilon^2\Delta v+v=|u|^{p-2}u &\text{in }A,\\
u,v>0\text{ in }A,\qquad u=v=0 &\text{on }\partial A,
\end{cases}
\end{equation}
where $A=\{x\in\mathbb R^{2m}:a<|x|<b\}$, $0<a<b$, $m\geq2$, and
$\varepsilon>0$ is a small parameter. Throughout the paper we assume
\begin{equation}\label{exponents}
2<p,q<\frac{2(m+1)}{m-1}.
\end{equation}
We seek positive solutions whose two components concentrate along the same
$(m-1)$-dimensional sphere in the inner boundary of $A$ as $\varepsilon\to0$.

We first recall some results for scalar equations. For the Dirichlet problem
\begin{equation}\label{1.2}
-\varepsilon^2\Delta u+u=f(u)\quad\text{in }D,\qquad u=0\quad\text{on }\partial D,
\end{equation}
Ni and Wei \cite{Ni-Wei} proved, under suitable assumptions on $f$, that ground states concentrate at points maximizing the distance to the boundary.
Del Pino and Felmer \cite{DelPino-Felmer} obtained related results for more general
nonlinearities. Concentration may also occur on sets of positive dimension; see
\cite{Ambrosetti-Malchiodi-Ni-I,Ambrosetti-Malchiodi-Ni,Esposito-Mancini-Santra-Srikanth}. For the connection
with reaction-diffusion models, we refer to \cite{Ni-WEI-2}.

For the radial problem
\begin{equation}\label{1.3}
-\varepsilon^2\Delta u+V(|x|)u=|u|^{p-2}u\quad\text{in }D,
\qquad u=0\quad\text{on }\partial D,
\end{equation}
with $D=\{x\in\mathbb R^N:a<|x|<1\}$ and a positive radial potential $V$,
the location of a radial layer is governed by
$$
M(r)=r^{N-1}V(r)^{\theta},\qquad
\theta=\frac{p}{p-2}-\frac12.
$$
Under the hypotheses of \cite{Ambrosetti-Malchiodi-Ni}, the condition $M'(a)>0$
leads to radial solutions concentrating on spheres whose radii tend to $a$.

Lower dimensional concentration sets can be obtained by exploiting the geometry
of the annulus. Ruf and Srikanth \cite{Ruf-Srikanth} used the Hopf fibration to
construct solutions concentrating on circles in dimension four. Their later work
\cite{Ruf-Srikanth-2014} treats concentration on Hopf fibers of dimensions one,
three and seven. Pacella and Srikanth \cite{Pacella-Srikanth} introduced a reduction
from dimension $2m$ to dimension $m+1$ and obtained positive and sign-changing
solutions concentrating on one or two $(m-1)$-dimensional spheres. Zhang and
do \'{O} \cite{zhang-doo} applied a related reduction to coupled Schr\"odinger
systems with cubic nonlinearities in a four-dimensional annulus.

For Hamiltonian systems, one must also deal with the strongly indefinite
variational structure. The energy associated with
\begin{equation}\label{1.6}
\begin{cases}
-\varepsilon^2\Delta u+u=g(v) &\text{in }D,\\
-\varepsilon^2\Delta v+v=f(u) &\text{in }D,\\
u=v=0 &\text{on }\partial D
\end{cases}
\end{equation}
is
$$
\int_D\bigl(\varepsilon^2\nabla u\cdot\nabla v+uv-F(u)-G(v)\bigr)\,\mathrm{d}x,
\qquad F(s)=\int_0^s f(t)\,\mathrm{d}t,\quad G(s)=\int_0^s g(t)\,\mathrm{d}t.
$$
Its quadratic part is positive on pairs $(u,u)$ and negative on pairs $(u,-u)$,
with $u\neq0$. Pistoia and Ramos \cite{PR} established point
concentration for least energy solutions of Dirichlet Hamiltonian systems with
subcritical nonlinearities. Manna and Sahoo \cite{Manna-Sahoo} studied Neumann
systems with positive variable coefficients and used dimensional reduction to
obtain higher dimensional concentration. On compact Riemannian manifolds without
boundary, Kannoth and Manna \cite{Kannoth-Manna} obtained point concentration at
maxima of the scalar curvature.

Concentration along positive dimensional sets also occurs under perturbations
of the exponents. Guo, Liu and Peng \cite{Guo-Liu-Peng} constructed boundary
layers for supercritical Lane-Emden systems on suitable symmetric domains.
Their solutions concentrate along products of spheres as the exponents approach
a critical hyperbola in a lower dimension. Chen and Wang \cite{Chen-Wang}
constructed solutions concentrating along minimal submanifolds of warped
product manifolds for asymptotically critical Hamiltonian systems.

More directly related to our problem is the recent preprint of Manna, Kannoth
and Ruf \cite{Manna-Kannoth-Ruf}. They consider Dirichlet Hamiltonian systems
with power weights on annuli in $\mathbb R^{2m}$, for $m=2,4,8$, and obtain
concentration on $(m-1)$-dimensional spheres. In particular, their Corollary~1.1
treats the constant-coefficient system \eqref{1.1} and gives concentration of
both components on a common sphere in the inner boundary. Their proof uses the
Hopf fibrations $\mathbb S^{2m-1}\to\mathbb S^m$ to lift solutions of a weighted
system in $\mathbb R^{m+1}$, without requiring axial symmetry of these solutions.
This construction accounts for the restriction $m=2,4,8$.
With the nonlinearities written as in \eqref{1.1}, their exponent condition reads
$$
p,q>2,\qquad \frac1p+\frac1q>\frac{m-1}{m+1}.
$$

Here we use the reduction of Pacella and Srikanth \cite{Pacella-Srikanth}, which
applies to every integer $m\geq2$ in the class of $O(m)\times O(m)$-invariant
functions. To use this correspondence for the system, we prove that the two
components of the reduced ground state are axially symmetric about the same
axis. This allows us to extend the sphere concentration result for \eqref{1.1}
to all $m\geq2$ under \eqref{exponents}. Our exponent range is narrower than
the preceding hyperbola condition, but still includes pairs that are critical
or supercritical in dimension $2m$, namely,
$$
\frac1p+\frac1q\leq\frac{m-1}{m}.
$$
In this regime, the Pohozaev identity rules out positive classical Dirichlet
solutions on smooth star-shaped domains; see \cite{Vandervorst}. The annular
geometry allows us to obtain the following result.

\begin{theo}\label{th1.1}

Let $2<p,q<2(m+1)/(m-1)$. Then, for all sufficiently small $\varepsilon>0$,
system \eqref{1.1} admits a positive solution $(u_\varepsilon,v_\varepsilon)$
which concentrates around an $(m-1)$-dimensional sphere on the inner boundary of $A$.
More precisely, writing $x=(y_1,y_2)\in\mathbb R^m\times\mathbb R^m$, the family
can be chosen invariant under separate rotations of $y_1$ and $y_2$. There are
radii $r_\varepsilon^u,r_\varepsilon^v\in(a,b)$ such that
\begin{equation}\label{radii}
r_\varepsilon^u,r_\varepsilon^v\to a,\qquad
\frac{r_\varepsilon^u-a}{\varepsilon},\,
\frac{r_\varepsilon^v-a}{\varepsilon}\to\infty,\qquad
\frac{|r_\varepsilon^u-r_\varepsilon^v|}{\varepsilon}\to0.
\end{equation}
The functions $u_\varepsilon$ and $v_\varepsilon$ attain their maxima on
$$
\Gamma_\varepsilon^u=\{(y_1,0):|y_1|=r_\varepsilon^u\},\qquad
\Gamma_\varepsilon^v=\{(y_1,0):|y_1|=r_\varepsilon^v\},
$$
respectively. Moreover, there are constants
$C,c>0$, independent of small $\varepsilon$, such that
\begin{equation}\label{liftdecay}
\begin{split}
u_\varepsilon(x)\leq C\exp\bigl(-c\,\operatorname{dist}(x,\Gamma_\varepsilon^u)/\varepsilon\bigr),\
v_\varepsilon(x)\leq C\exp\bigl(-c\,\operatorname{dist}(x,\Gamma_\varepsilon^v)/\varepsilon\bigr).
\end{split}
\end{equation}
\end{theo}

\begin{Remark}\label{r1.1}

Condition \eqref{exponents} is subcritical in dimension $m+1$, but allows critical
and supercritical pairs in dimension $2m$. For example, we may take $(p,q)=(4,5)$
when $m=2$ and $(p,q)=(3,10/3)$ when $m=3$.
The concentration sphere lies on the limiting boundary, whereas the maxima remain
inside $A$ and their distance to the boundary is large compared with $\varepsilon$.
Thus \eqref{radii} is consistent with the Dirichlet condition.

When $p=q=s$, subtracting the two equations in \eqref{1.1} and testing with
$u-v$ gives
$$
\begin{aligned}
\varepsilon^2\int_A|\nabla(u-v)|^2\,\mathrm{d}x
+\int_A|u-v|^2\,\mathrm{d}x
&=-\int_A(u^{s-1}-v^{s-1})(u-v)\,\mathrm{d}x\leq0.
\end{aligned}
$$
Hence $u=v$. The existence of sphere-concentrating solutions in this case
follows from the scalar result of Pacella and Srikanth \cite{Pacella-Srikanth}.
The system argument below also allows unequal exponents.

\end{Remark}

We now outline the proof. In the class of functions invariant under separate
rotations of $y_1$ and $y_2$, the reduction of \cite{Pacella-Srikanth} transforms
\eqref{1.1} into the system
\begin{equation}\label{2.3}
\begin{cases}
-\varepsilon^2\Delta u+\frac{1}{2|x|}u=\frac{1}{2|x|}|v|^{q-2}v &\text{in }\Omega,\\[3pt]
-\varepsilon^2\Delta v+\frac{1}{2|x|}v=\frac{1}{2|x|}|u|^{p-2}u &\text{in }\Omega,\\[3pt]
u,v>0\text{ in }\Omega,\qquad u=v=0 &\text{on }\partial\Omega,
\end{cases}
\end{equation}
where $\Omega=\{x\in\mathbb R^{m+1}:a^2/2<|x|<b^2/2\}$.

To lift a solution of \eqref{2.3}, both components must have the same axis of
symmetry. Direct polarization of the indefinite quadratic form does not provide
the required comparison. We use the inversion formulation of Bonheure,
Moreira dos Santos and Ramos \cite{Bonheure-Ramos-ground} and adapt their
polarization argument in \cite{Bonheure-Ramos-symmetry} to the weighted annular
problem; see also the survey \cite{Bonheure-Tavares}. Ground states are
characterized by a best constant, and polarization of the inverse Dirichlet
operator gives foliated Schwarz symmetry with a common direction.

To study concentration, we first prove a quadratic error estimate for the maximum
of the energy on a generalized Nehari fiber, with constants independent of the
rescaled domain. Applied to cut-off ground states of the frozen-coefficient
problem, this estimate gives the upper energy bound. We next use a half-space
nonexistence result to exclude boundary profiles. Comparison with the lower
energy bound forces concentration near the inner boundary and gives strong
convergence of the rescaled solutions. A local estimate for their zero extensions
then yields uniform exponential decay. Finally, strict radial monotonicity of
the limiting profiles shows that the maxima of the two components are at
distance $o(\varepsilon)$ from each other.

For concentration under Neumann boundary conditions, we refer to
\cite{Avila-Yang,Ramos-Yang,Ni-Takagi-1993,Guo2025JDE}; for related problems in the
whole space, see
\cite{Rabinowitz1992,wangxuefeng1993,Ambrosetti-Malchiodi-2007,Byeon-Jeanjean-2007}.

The paper is organized as follows. Section \ref{sec:2} recalls the dimensional
reduction. Section \ref{sec:3} establishes the inversion formulation and the
symmetry of ground states. Section \ref{sec:4} gives the generalized Nehari
characterization and the uniform fiber estimate. The upper energy bound is proved
in Section \ref{sec:5}. Section \ref{sec:6} contains the concentration analysis and
the proof of Theorem~\ref{th1.1}.

\vskip 0.2in
\noindent\textbf{Notation.}
From Section \ref{sec:2} onward,
\begin{itemize}
\item $N=m+1$, $2^*=2N/(N-2)$ and $w(x)=(2|x|)^{-1}$ on
$\overline\Omega$. In particular, $b^{-2}\leq w\leq a^{-2}$.
\item We use the norm $\|u\|_{H^1(D)}^2=\int_D(|\nabla u|^2+u^2)\mathrm{d}x$ and the corresponding inner product in $H^1(D)$. And we use the norm $\|u\|_{L^s(D)}=\bigl(\int_D|u|^sdx\bigr)^{\frac1s}$ in $L^s(D)$ and for briefly  we denote $\|u\|_{s}$ as $\|u\|_{L^s(\mathbb{R}^N)}$.
\item Functions in $H^1_0(D)$ are extended by zero when compared
on $\mathbb R^N$.
\item The notation $B_R(y)$ stands for the open ball of radius $R$
centered at $y$.
\item The letter $C$ denotes a positive constant whose value may change
from line to line. Unless otherwise specified, constants in uniform estimates
are independent of sufficiently small $\varepsilon$.
\end{itemize}

\section{Axial symmetry and dimensional reduction}\label{sec:2}

We recall the reduction of \cite{Pacella-Srikanth} and specify the symmetry class
in which it is used. Let $D$ be a ball or an annulus centered at the origin in
$\mathbb R^N$.

\begin{Definition}\label{def2.1}
A function $f:D\to\mathbb R$ is axially symmetric with respect to a unit vector
$e$ if it depends only on $\rho=|z|$ and
$\varphi=\arccos(z\cdot e/|z|)$.
\end{Definition}
\begin{Definition}\label{def2.2}
An axially symmetric function is foliated Schwarz symmetric with respect to $e$
if it is non-increasing in $\varphi$.
\end{Definition}

The value at the origin, when $0\in D$, is understood by continuity. After a
rotation we may take $e=e_N$. For an axially symmetric $C^2$ function,

$$
\Delta_{\mathbb R^N}f=f_{\rho\rho}+\frac{N-1}{\rho}f_\rho+
\frac1{\rho^2}\bigl(f_{\varphi\varphi}+(N-2)\cot\varphi\,f_\varphi\bigr).
$$

Write $x=(y_1,y_2)\in\mathbb R^m\times\mathbb R^m$ and set
$$
|y_1|=r\cos\theta,\qquad |y_2|=r\sin\theta,\qquad
r=|x|,\quad 0\leq\theta\leq\pi/2.
$$
A function invariant under separate rotations of $y_1$ and $y_2$ depends only on
$(r,\theta)$. Its Laplacian is
$$
\Delta_{\mathbb R^{2m}}f=f_{rr}+\frac{2m-1}{r}f_r+
\frac1{r^2}\bigl(f_{\theta\theta}+(m-1)(\cot\theta-\tan\theta)f_\theta\bigr).
$$
For such a function define $\widetilde f$ by

\begin{equation}\label{2.7}
\rho=\frac{r^2}{2},\qquad\varphi=2\theta,\qquad
\widetilde f(\rho,\varphi)=f(\sqrt{2\rho},\varphi/2).
\end{equation}

A direct computation gives
\begin{equation}\label{lapreduction}
\Delta_{\mathbb R^{2m}}f
=2\rho\left[\widetilde f_{\rho\rho}+\frac m\rho\widetilde f_\rho+
\frac1{\rho^2}\bigl(\widetilde f_{\varphi\varphi}+(m-1)\cot\varphi\,
\widetilde f_\varphi\bigr)\right]
=2\rho\,\Delta_{\mathbb R^{m+1}}\widetilde f.
\end{equation}
The computation is valid for every integer $m\geq2$.
Consequently, \eqref{1.1} in this symmetry class is equivalent to \eqref{2.3}
in the class of pairs axially symmetric with respect to $e_N$. Conversely, a
classical solution of \eqref{2.3} with this common axis lifts by \eqref{2.7} to a
solution of \eqref{1.1}. The regularity at $\theta=0,\pi/2$ follows from the axial
regularity at $\varphi=0,\pi$; this is the correspondence established in
\cite{Pacella-Srikanth}. In particular, a point $\rho e_N$ corresponds to the sphere
\begin{equation}\label{spherecorrespondence}
\{(y_1,0):|y_1|=\sqrt{2\rho}\}.
\end{equation}
We shall apply this correspondence after rotating the common symmetry direction
of each reduced solution to $e_N$.

\section{Reduction by inversion and symmetry of ground states}\label{sec:3}

Throughout this section $\varepsilon>0$ is fixed. We retain it in the operator
$$
L_\varepsilon=-\varepsilon^2\Delta+w(x),\qquad q'=\frac q{q-1},\qquad
\delta_{p,q}=1-\frac1p-\frac1q>0.
$$
For a smooth bounded connected domain $D\subseteq\Omega$, set
$X_D=W^{2,q'}(D)\cap W^{1,q'}_0(D)$. We first record the estimates needed for
the inversion formulation.

\begin{Lemma}[Dirichlet estimates]\label{dirichletestimate}
For every $u\in X_D$, we have
\begin{equation}\label{dirichletbounds}
\begin{split}
\|u\|_{L^{q'}(D)}&\leq b^2\|L_\varepsilon u\|_{L^{q'}(D)},\\
\|u\|_{W^{2,q'}(D)}&\leq C_{\varepsilon,D}
\|L_\varepsilon u\|_{L^{q'}(D)}.
\end{split}
\end{equation}
The first constant is independent of $\varepsilon$ and $D\subseteq\Omega$.
The second constant may depend on $\varepsilon$ and $D$.
\end{Lemma}
\begin{proof}
 Since $q<2^*$, we have
$q'>2N/(N+2)$. Sobolev embedding and the zero trace therefore give
$u\in H^1_0(D)$, and $L_\varepsilon u\in L^{q'}(D)\subset H^{-1}(D)$.
To justify testing with $|u|^{q'-2}u$, for $\eta>0$ define
$$
\psi_\eta(t)=t(t^2+\eta^2)^{({q'}-2)/2},
$$
which is Lipschitz, vanishes at zero, and satisfies
$$
\psi_\eta'(t)=(t^2+\eta^2)^{({q'}-4)/2}
\bigl(\eta^2+({q'}-1)t^2\bigr)\geq0.
$$
Thus $\psi_\eta(u)\in H^1_0(D)$ is an admissible test function. We obtain
$$
\varepsilon^2\int_D\psi_\eta'(u)|\nabla u|^2\,\mathrm{d}x
+\int_D wu\psi_\eta(u)\,\mathrm{d}x
=\int_D L_\varepsilon u\psi_\eta(u)\,\mathrm{d}x.
$$
 Since
$|\psi_\eta(u)|\leq |u|^{q'-1}$, dominated convergence as $\eta\to0$, the lower
bound $w\geq b^{-2}$, and H\"older's inequality yield
$$
b^{-2}\|u\|_{L^{q'}(D)}^{q'}\leq\int_D w|u|^{q'}\,\mathrm{d}x
\leq\int_D L_\varepsilon u|u|^{{q'}-2}u\,\mathrm{d}x
\leq\|L_\varepsilon u\|_{_{L^{q'}(D)}}\|u\|_{_{L^{q'}(D)}}^{{q'}-1}.
$$
This proves the first estimate, including the case $u=0$.

Since $D$ is smooth and $u$ has zero Dirichlet trace, the global elliptic
estimate \cite{GT} gives
$$
\|u\|_{W^{2,{q'}}(D)}
\leq C_{\varepsilon,D}\bigl(\|L_\varepsilon u\|_{L^{q'}(D)}+\|u\|_{L^{q'}(D)}\bigr).
$$
Substituting $\|u\|_{L^{q'}(D)}\leq b^2\|L_\varepsilon u\|_{L^{q'}(D)}$ into this inequality and changing the
constant proves the second estimate in \eqref{dirichletbounds}.
Here $\varepsilon$ and $D$ are fixed;
no uniform bound for $C_{\varepsilon,D}$ is asserted.
\end{proof}

We now define
\begin{equation}\label{bestconstant}
\begin{split}
A_D(u)&=\int_D w^{-1/(q-1)}|L_\varepsilon u|^{q'}\,\mathrm{d}x,
\qquad B_D(u)=\int_D w|u|^p\,\mathrm{d}x,\\
S_{\varepsilon,D}&=\inf_{u\in X_D\setminus\{0\}}
\frac{A_D(u)}{B_D(u)^{q'/p}}.
\end{split}
\end{equation}
The corresponding system is 
\begin{equation}\label{bu}
\begin{cases}
-\varepsilon^2\Delta u+\frac{1}{2|x|}u=\frac{1}{2|x|}|v|^{q-2}v &\text{in }D,\\[3pt]
-\varepsilon^2\Delta v+\frac{1}{2|x|}v=\frac{1}{2|x|}|u|^{p-2}u &\text{in }D,\\[3pt]
u,v>0\text{ in }D,\qquad u=v=0 &\text{on }\partial D.
\end{cases}
\end{equation}
The energy of system (\ref{bu}) is
\begin{equation}\label{energy}
I_{\varepsilon,D}(u,v)=\int_D(\varepsilon^2\nabla u\cdot\nabla v+wuv)\,\mathrm{d}x
-\int_D w\left(\frac{|u|^p}{p}+\frac{|v|^q}{q}\right)\,\mathrm{d}x.
\end{equation}
A  nontrivial solution $z$  of (\ref{bu}) is a ground state if $I_{\varepsilon,D}(z)$ is the least energy among
all nontrivial critical points of (\ref{bu}), and denote this least energy by $c_{\varepsilon,D}$.
For $D=\Omega$ we write $I_\varepsilon$ and $c_\varepsilon$.

The next lemma adapts the inversion formulation in
\cite[Section 2]{Bonheure-Ramos-ground} to $L_\varepsilon$. We include the
proof, in particular to recover the Dirichlet trace of the second component.

\begin{Lemma}\label{groundinversion}
The constant $S_{\varepsilon,D}$ is positive and attained. Moreover,
\begin{equation}\label{energyconstant}
c_{\varepsilon,D}=\delta_{p,q}S_{\varepsilon,D}^{p/(p-q')}.
\end{equation}
System (\ref{bu}) has a ground state, and both components are
strictly positive. The first component of every ground state, normalized by
$B_D(u)=1$, is a minimizer of \eqref{bestconstant}.
\end{Lemma}
\begin{proof}
We divide the proof into three steps.

\noindent\textit{Step 1: attainment and the weak inversion formulation.}
Condition \eqref{exponents} implies
$1/p+1/q>(N-2)/N$, hence $X_D$ is compactly embedded into $L^p(D)$.
Lemma~\ref{dirichletestimate}, the bounds on $w$, and the continuous embedding
into $L^p(D)$ show that $S_{\varepsilon,D}>0$.
Minimizing $A_D$ under $B_D=1$ and using weak
lower semicontinuity and the compact embedding yields a minimizer $u$.
Put $S=S_{\varepsilon,D}$ and
$v=w^{-1/(q-1)}|L_\varepsilon u|^{q'-2}L_\varepsilon u$. The Euler equation is
\begin{equation}\label{weaknavier}
\int_D vL_\varepsilon\phi\,\mathrm{d}x
=S\int_D w|u|^{p-2}u\phi,\quad \phi\in X_D,
\end{equation}
and $L_\varepsilon u=w|v|^{q-2}v$ almost everywhere. The expression
$|s|^{q'-2}s$ is defined to be zero at $s=0$.

It is essential that \eqref{weaknavier} holds for all $\phi\in X_D$, rather than
only for compactly supported tests. Indeed, $u\in L^p(D)$ and $v\in L^q(D)$.
Since $p,q<2^*$, the nonlinear terms belong to $H^{-1}(D)$. Let
$\bar v\in H^1_0(D)$ solve $L_\varepsilon\bar v=Sw|u|^{p-2}u$.
For $\psi\in C_c^\infty(D)$, the Dirichlet solution of
$L_\varepsilon\phi=\psi$ belongs to $X_D\cap H^1_0(D)$. Substitution in
\eqref{weaknavier} gives $\int_D(v-\bar v)\psi\mathrm{d}x=0$. Thus $v=\bar v\in H^1_0(D)$.
The first component already belongs to $H^1_0(D)$ by the embedding used in
Lemma~\ref{dirichletestimate}. We have thus verified the Dirichlet condition
on the second component, which is the natural boundary
condition in the inversion formulation.

The pair solves
\begin{equation}\label{normalizedsystem}
L_\varepsilon u=w|v|^{q-2}v,\qquad
L_\varepsilon v=Sw|u|^{p-2}u,\qquad u=v=0\text{ on }\partial D.
\end{equation}
Subcritical elliptic bootstrap gives boundedness and then classical regularity.
For completeness, one may apply the scalar subsolution estimate to $|u|+|v|$:
its nonlinear coefficient is bounded by $C(1+(|u|+|v|)^{s-2})$, where
$s=\max\{p,q\}<2^*$. This coefficient belongs locally to $L^t$ for some $t>N/2$,
since $u,v\in H^1_0(D)$. Interior and Dirichlet boundary estimates then yield
$u,v\in L^\infty(D)$ and $C^{2,\alpha}(\overline D)$ for some $\alpha\in(0,1)$.

\noindent\textit{Step 2: positivity.}
Let $\widehat{u}$ solve $L_\varepsilon\widehat{u}=|L_\varepsilon u|$ with zero
boundary values. Since
$L_\varepsilon(\widehat u\pm u)=|L_\varepsilon u|\pm L_\varepsilon u\geq0$ in $D$
and $\widehat u\pm u=0$ on $\partial D$, the maximum principle gives
$\widehat u\geq|u|$. Moreover, $A_D(\widehat u)=A_D(u)=S$ and
$B_D(\widehat u)\geq B_D(u)=1$. The minimality of the quotient forces
$B_D(\widehat u)=1$, and hence $\widehat u=|u|$.
If $L_\varepsilon u$ had both a nontrivial positive part and a nontrivial negative part,
the strong maximum principle applied to
$\widehat u+u$ and $\widehat u-u$ would give $\widehat u>|u|$ in $D$, a
contradiction. Hence, after changing the sign, $L_\varepsilon u\geq0$ and
$u>0$. It follows that $v\geq0$; the second equation in
\eqref{normalizedsystem} then gives $v>0$ and $L_\varepsilon u>0$.

\noindent\textit{Step 3: the energy identity.}
 For any nontrivial solution $(U,V)$ of (\ref{bu}), testing its
equations with $V$ and $U$, respectively, gives
\begin{equation}\label{commonmass}
K:=\int_D(\varepsilon^2\nabla U\cdot\nabla V+wUV)\mathrm{d}x
=\int_Dw|U|^p\mathrm{d}x=\int_Dw|V|^q\mathrm{d}x>0.
\end{equation}
Thus $I_{\varepsilon,D}(U,V)=\delta_{p,q}K$ and
$S\leq K^{1-q'/p}$. So
$c_{\varepsilon,D}\geq\delta_{p,q}S^{p/(p-q')}$. To attain this value, take a minimizer $u$ from Step 1
and set $v=w^{-1/(q-1)}|L_\varepsilon u|^{q'-2}L_\varepsilon u$.
By Step 2, we may choose $u,v>0$. Note that $(u,v)$ solves \eqref{normalizedsystem}, we know
$$
U_0=tu,\qquad V_0=t^{1/(q-1)}v,\qquad t=S^{1/(p-q')}
$$
solves (\ref{bu}).
Since $K=t^p=S^{p/(p-q')}$, we get $c_{\varepsilon,D}\leq I_{\varepsilon,D}(U_0,V_0)=\delta_{p,q}S^{p/(p-q')}$. Therefore, (\ref{energyconstant}) yields. This ends the proof.
\end{proof}

We now take $D=\Omega$. Let $\mathcal{H}$ be the family of closed half-spaces
whose boundary contains the origin. For $H\in\mathcal H$, denote reflection
across $\partial H$ by $\sigma_H$, and write $\bar x=\sigma_H(x)$ and
$\bar f(x)=f(\bar x)$. The polarization of $f$ is

\begin{equation}\label{3.0.0}
f_H(x)=\begin{cases}
\max\{f(x),f(\bar x)\},&x\in H\cap\Omega,\\
\min\{f(x),f(\bar x)\},&x\in\Omega\setminus H.
\end{cases}
\end{equation}

The restriction $0\in\partial H$ ensures that both $\Omega$ and $w$ are invariant
under $\sigma_H$. The comparison below is the version for $L_\varepsilon$ of
\cite[Lemma 3.7]{Bonheure-Ramos-symmetry}; the underlying polarization
comparison goes back to Brock and Solynin \cite[Section 9]{Brock-Solynin}.

\begin{Lemma}\label{inversepolar}
Let $f\in L^\infty(\Omega)$ and $H\in\mathcal H$. Let $u,v$ be the Dirichlet
solutions of $L_\varepsilon u=f$ and $L_\varepsilon v=f_H$. Then
$$
v=v_H,\qquad v+\bar v=u+\bar u\quad\text{in }\Omega,
\qquad v\geq u_H\quad\text{in }H\cap\Omega.
$$
For every $\phi\in L^\infty(\Omega)$,
\begin{equation}\label{linearpolar}
\int_\Omega u\phi \mathrm{d}x\leq\int_\Omega v\phi_H\mathrm{d}x.
\end{equation}
If $f\geq0$, $\phi\geq0$ and $s>1$, then
\begin{equation}\label{convexpolar}
\int_\Omega u^s\phi\mathrm{d}x\leq\int_\Omega v^s\phi_H\mathrm{d}x.
\end{equation}
\end{Lemma}
\begin{proof}
The maximum principle applied to $v-\bar v$ on $H\cap\Omega$ gives
$v\geq\bar v$. Reflection interchanges $(\operatorname{int}H)\cap\Omega$ and $\Omega\setminus H$ and fixes $\partial H\cap\Omega$. Hence, $v=v_H$ in $\Omega$. Since $f+\bar f=f_H+\overline{f_H}$, uniqueness of the
Dirichlet problem gives $v+\bar v=u+\bar u$. In particular $v=u=\bar u$ on
$\partial H\cap\Omega$, and thus $v-u=0, v-\bar{u}=0$ on $\partial(H\cap \Omega)$.
Applying the maximum principle to $v-u$ and $v-\bar u$ on $H\cap\Omega$ gives $v\geq\max\{u,\bar u\}$ there. So $v\geq u_H$ in $H\cap \Omega$.
We note that
$$\int_\Omega u\phi\mathrm{d}x-\int_\Omega v\phi_H\mathrm{d}x=\int_{H\cap \Omega} [u\phi+\bar u\bar\phi-v\phi_H-\bar v\overline{\phi_H}]\mathrm{d}x.$$
Then by the identity
$$
u\phi+\bar u\bar\phi-v\phi_H-\bar v\overline{\phi_H}
=(\phi-\phi_H)(v-\bar u)+(\bar\phi-\phi_H)(v-u)\leq0
$$
on $H\cap\Omega$, we obtain \eqref{linearpolar}.

Suppose now that $f,\phi\geq0$. At each reflected pair set
$a=\max\{u,\bar u\}$, $b=\min\{u,\bar u\}$, $A=v$, $B=\bar v$ on
$H\cap\Omega$. Then $A\geq a\geq b\geq B\geq0$ and $A+B=a+b$.
If $\eta=\max\{\phi,\bar\phi\}$ and $\xi=\min\{\phi,\bar\phi\}$, the two-term
rearrangement inequality gives
$$u^s\phi+\bar u^s\bar\phi\leq\eta a^s+\xi b^s.$$
Let $\delta:=A-a=b-B\geq 0$ and introduce
$$F(t):=(a+t)^s+(b-t)^s,\qquad 0\leq t\leq \delta.$$
Since $0\leq\delta\leq b$, $s>1$ and $a+t\geq b-t\geq0$, we have
$$F'(t)=s[(a+t)^{s-1}-(b-t)^{s-1}]\geq0.$$
Thus $F(\delta)\geq F(0)$, that is,
$$A^s+B^s\geq a^s+b^s.$$
Using $A\geq a$ and $\eta\geq\xi\geq0$, we obtain
$$
\begin{aligned}
\eta a^s+\xi b^s
&=(\eta-\xi)a^s+\xi(a^s+b^s)\\
&\leq(\eta-\xi)A^s+\xi(A^s+B^s)
=\eta A^s+\xi B^s.
\end{aligned}
$$
Finally, the identity
$$
\begin{aligned}
&\int_\Omega u^s\phi\,\mathrm{d}x-\int_\Omega v^s\phi_H\,\mathrm{d}x=\int_{H\cap\Omega}
\big[u^s\phi+\bar u^s\bar\phi-v^s\phi_H
-\bar v^s\overline{\phi_H}\big]\,\mathrm{d}x
\end{aligned}
$$
yields \eqref{convexpolar}.
\end{proof}

As in \cite[Lemmas 3.4--3.5]{Bonheure-Ramos-symmetry}, define
$T:X_\Omega\to X_\Omega^*$ by
$$
\langle T(u),\phi\rangle
=\int_\Omega w^{-1/(q-1)}|L_\varepsilon u|^{q'-2}
L_\varepsilon u\,L_\varepsilon\phi\,\mathrm{d}x.
$$
Let $\mathcal L:X_\Omega\to L^{q'}(\Omega)$ be the Dirichlet isomorphism
$\mathcal Lu=L_\varepsilon u$. Its Banach adjoint
$\mathcal L':L^q(\Omega)\to X_\Omega^*$ is also an isomorphism and satisfies
$$
\langle\mathcal L'v,\phi\rangle
=\int_\Omega vL_\varepsilon\phi\,\mathrm{d}x,
\qquad v\in L^q(\Omega),\quad\phi\in X_\Omega.
$$
The weighted power map
$$
\mathcal N:L^{q'}(\Omega)\to L^q(\Omega),\qquad
\mathcal N(f)=w^{-1/(q-1)}|f|^{q'-2}f,
$$
has the inverse $\mathcal N^{-1}(v)=w|v|^{q-2}v$.
Since $w$ is bounded above and away from zero, both maps are continuous
between the indicated spaces. Thus $\mathcal N$ is a homeomorphism, and
$$
T=\mathcal L'\circ\mathcal N\circ\mathcal L,
\qquad
T^{-1}=\mathcal L^{-1}\circ\mathcal N^{-1}\circ(\mathcal L')^{-1}.
$$
In particular, $T:X_\Omega\to X_\Omega^*$ is a homeomorphism.

For $h\in L^\infty(\Omega)$, we identify $h$ with the functional
$\phi\mapsto\int_\Omega h\phi\,\mathrm{d}x$ on $X_\Omega$.
To identify $T^{-1}h$, solve successively
$$
\begin{cases}
L_\varepsilon v=h&\text{in }\Omega,\\
v=0&\text{on }\partial\Omega,
\end{cases}
\qquad
\begin{cases}
L_\varepsilon u=w|v|^{q-2}v&\text{in }\Omega,\\
u=0&\text{on }\partial\Omega.
\end{cases}
$$
Here $v\in W^{2,q}(\Omega)\cap W^{1,q}_0(\Omega)$ and $u\in X_\Omega$.
The differential expression $L_\varepsilon=-\varepsilon^2\Delta+w$ is
formally self-adjoint. Green's formula and the zero Dirichlet traces therefore
give, for every $\phi\in X_\Omega$,
$$
\begin{aligned}
\int_\Omega vL_\varepsilon\phi\,\mathrm{d}x
&=\varepsilon^2\int_\Omega\nabla v\cdot\nabla\phi\,\mathrm{d}x
+\int_\Omega wv\phi\,\mathrm{d}x\\
&=\int_\Omega(L_\varepsilon v)\phi\,\mathrm{d}x
=\int_\Omega h\phi\,\mathrm{d}x.
\end{aligned}
$$
Since $\mathcal N(L_\varepsilon u)=v$, it follows that
$$
\langle T(u),\phi\rangle
=\int_\Omega vL_\varepsilon\phi\,\mathrm{d}x
=\int_\Omega h\phi\,\mathrm{d}x.
$$
Hence $T(u)=h$, and injectivity of $T$ shows that the two Dirichlet problems
give its unique solution. This construction also identifies
$\mathcal N(L_\varepsilon u)$ with the zero Dirichlet solution $v$, providing
its regularity and boundary trace.

If $h\geq0$ and $T(\widetilde u)=h_H$, Lemma~\ref{inversepolar}, with $s=q$
and $\phi=w$ in \eqref{convexpolar}, gives
\begin{equation}\label{nonlinearpolar}
\langle T(u),u\rangle=\int_\Omega wv^q\mathrm{d}x
\leq\int_\Omega w\widetilde v^q\mathrm{d}x
=\langle T(\widetilde u),\widetilde u\rangle.
\end{equation}

\begin{Proposition}\label{commonsymmetry}
Every positive ground state of \eqref{2.3} is foliated Schwarz symmetric with
respect to a common unit vector $e\in\mathbb R^N$.
\end{Proposition}
\begin{proof}
We adapt the proof of \cite[Theorem 1.3]{Bonheure-Ramos-symmetry}.
By Lemma~\ref{groundinversion}, it is enough to consider a positive minimizer $u$
of $S=S_{\varepsilon,\Omega}$ with $B_\Omega(u)=1$, and its companion
$v=w^{-1/(q-1)}(L_\varepsilon u)^{1/(q-1)}$. They satisfy
\eqref{normalizedsystem}. Positive scalar multipliers do not affect the symmetry.

\noindent\textit{Step 1: polarization preserves the minimizing property.}
For $H\in\mathcal H$, let $T(\widetilde u)=Swu_H^{p-1}$. Since $w_H=w$,
\eqref{nonlinearpolar} and H\"older's inequality yield
$$
S\leq\langle T(\widetilde u),\widetilde u\rangle
=S\int_\Omega w\widetilde u\,u_H^{p-1}\mathrm{d}x
\leq S B_\Omega(\widetilde u)^{1/p},
$$
where we have used the facts $\widetilde u\geq0$ and $B_\Omega(u_H)=1$. Consequently
$B_\Omega(\widetilde u)\geq1$, whereas the definition of $S$ gives
$$
S\leq\frac{\langle T(\widetilde u),\widetilde u\rangle}
{B_\Omega(\widetilde u)^{q'/p}}
\leq S B_\Omega(\widetilde u)^{(1-q')/p}\leq S.
$$
All inequalities are equalities. Hence $B_\Omega(\widetilde u)=1$,
$\widetilde u$ is a minimizer, and its Euler equation implies
$Sw\widetilde u^{p-1}=T(\widetilde u)=Swu_H^{p-1}$. Thus
$\widetilde u=u_H$. In particular, $u_H\in X_\Omega$ and is a minimizer.

\noindent\textit{Step 2: the same direction works for both components.}
Fix $r\in(a^2/2,b^2/2)$ and choose $y$ maximizing $u$ on $|x|=r$.
Set $e=y/r$. Let $H\in\mathcal H$ contain $e$ in its interior and let
$v^H=w^{-1/(q-1)}(L_\varepsilon u_H)^{1/(q-1)}$. By Step 1,
$$
L_\varepsilon u_H=w(v^H)^{q-1},\qquad
L_\varepsilon v^H=Swu_H^{p-1}.
$$
Put $f:=Sw u^{p-1}$, then we note that $f_H=Swu_H^{p-1}$.
So, by Lemma~\ref{inversepolar},
 $$v^H=(L_\varepsilon)^{-1} f_H \geq \big((L_\varepsilon)^{-1}f\big)_H=v_H~\hbox{on}~H\cap\Omega.$$
Since $v_H=\max\{v,\bar{v}\}\geq v$, we obtain that
$$v^H\geq v_H\geq v~ \hbox{on}~ H\cap\Omega.$$
 Therefore
$$
u_H-u\geq0,\qquad L_\varepsilon(u_H-u)\geq0
\quad\text{on }H\cap\Omega.
$$
Since $y\in H\cap \Omega$, $u_H(y)=\max\{u(y), u(\bar{y})\}=u(y)$ and $H\cap\Omega$ is connected, the strong maximum principle gives $u_H=u$ there. Their equations then imply $v^H=v$, so $v_H=v$ on $H\cap\Omega$ as well. The reflected identities give equality
throughout $\Omega$. If $e\in\partial H$, approximate $H$ by halfspaces
$H_j\in\mathcal H$ containing $e$ in their interiors, with the unit normals
converging to that of $H$. Continuity of $u,v$ and of the reflections then gives
$u_H=u$ and $v_H=v$ also in this case.

Thus $u_H=u$ and $v_H=v$ for every $H\in\mathcal H$ containing $e$.
We include the short geometric argument yielding foliated Schwarz symmetry
from these identities; see also \cite{Bartsch-Weth-Willem}.
Let $g$ denote either $u$ or $v$, and take $x,z\in\Omega$ with
$|x|=|z|$ and $x\cdot e\geq z\cdot e$. If $x\neq z$, set
$$
n=\frac{x-z}{|x-z|},\qquad
H=\{\zeta\in\mathbb R^N:\zeta\cdot n\geq0\}.
$$
Then $e\in H$. Since $|x|=|z|$, we have
$$
x\cdot n=\frac{|x-z|}{2}>0,
\qquad \sigma_Hx=x-2(x\cdot n)n=z.
$$
Consequently,
$$
g(x)=g_H(x)=\max\{g(x),g(z)\}\geq g(z).
$$
The same inequality is trivial when $x=z$. If $x\cdot e=z\cdot e$,
interchanging $x$ and $z$ gives $g(x)=g(z)$. Hence $g$ depends only on
$r=|x|$ and $\theta=\arccos(x\cdot e/|x|)$ and is non-increasing in
$\theta$ for each fixed $r$. Both $u$ and $v$ are therefore foliated Schwarz
symmetric with respect to the same unit vector $e$.
\end{proof}

\section{The generalized Nehari characterization and a uniform estimate}\label{sec:4}

Let $E_D=H^1_0(D)\times H^1_0(D)$, with
$$
\|(u,v)\|_{E_D}^2=\|u\|_{H^1(D)}^2+\|v\|_{H^1(D)}^2.
$$
We use the decomposition $E_D=E_D^+\oplus E_D^-$, where
$E_D^+=\{(\phi,\phi)\}$ and $E_D^-=\{(\phi,-\phi)\}$, with
$\phi\in H^1_0(D)$. For $z\notin E_D^-$, its fiber is
$$
\widehat E_D(z)=\{tz+(\phi,-\phi):t\geq0,\ \phi\in H^1_0(D)\}.
$$
For the energy $I_\varepsilon$ in \eqref{energy}, define
$$
\mathcal N_\varepsilon=\{z\in E_\Omega\setminus E_\Omega^-:
\langle I_\varepsilon'(z),z\rangle=0,\quad
\langle I_\varepsilon'(z),h\rangle=0\text{ for all }h\in E_\Omega^-\}.
$$
We first record the scalar inequality which will also give a uniform estimate
for almost critical points. For $s>2$, put
\begin{equation}\label{Hs}
H_s(t)=\frac{t^2+1}{2}-\frac1s-\frac{s-1}{s}t^{s/(s-1)},\qquad t\geq0.
\end{equation}
Since $H_s'(t)=t-t^{1/(s-1)}$, the function $H_s$ is strictly
decreasing on $(0,1)$ and strictly increasing on $(1,\infty)$, with $H_s(1)=0$.
Moreover,
$$
\lim_{t\to1}\frac{H_s(t)}{(t-1)^2}=\frac{s-2}{2(s-1)},\qquad
\lim_{t\to\infty}\frac{H_s(t)}{(t-1)^2}=\frac12,\qquad
H_s(0)=\frac12-\frac1s>0.
$$
The quotient therefore extends to a positive continuous function at $t=1$
and is bounded away from zero on $[0,\infty)$. Hence there is $\kappa_s>0$ such that
\begin{equation}\label{Hlower}
H_s(t)\geq\kappa_s(t-1)^2\qquad(t\geq0).
\end{equation}

For an open set $D\subseteq\mathbb R^N$, consider the functional
\begin{equation}\label{generalPhi}
\Phi(u,v)=\int_D(\nabla u\cdot\nabla v+A(x)uv)\,\mathrm{d}x
-\int_D B(x)\left(\frac{|u|^p}{p}+\frac{|v|^q}{q}\right)\,\mathrm{d}x,
\end{equation}
where $A,B\in L^\infty(D)$ satisfy $0<\kappa\leq A,B\leq K$ almost everywhere,
and $p,q$ satisfy \eqref{exponents}. Write
$$
(f,g)_A=\int_D(\nabla f\cdot\nabla g+Afg)\,\mathrm{d}x,
\qquad \|\phi\|_A^2=(\phi,\phi)_A.
$$
The basic fiber comparison is given in \cite[Lemma 3.3]{Tang-fiber}.
For the pure powers in \eqref{generalPhi}, we retain the remainder terms
involving $H_p$ and $H_q$. Their quadratic lower bounds will be used for almost
critical points on varying domains.

\begin{Lemma}[A fiber inequality]\label{fiberestimate}
Let $\Phi$ be given by \eqref{generalPhi}. For every $z=(u,v)\in E_D$,
$\phi\in H^1_0(D)$ and $t\geq0$, setting $h=(\phi,-\phi)$, we have
\begin{equation}\label{fiberineq}
\begin{split}
\Phi(tz+h)-\Phi(z)
\leq{}&\left\langle\Phi'(z),\frac{t^2-1}{2}z+th\right\rangle
-\|\phi\|_A^2\\
&-H_p(t)\int_D B|u|^p\,\mathrm{d}x
-H_q(t)\int_D B|v|^q\,\mathrm{d}x.
\end{split}
\end{equation}
The same estimate holds for $I_\varepsilon$, with $D=\Omega$, $B=w$ and
$\|\phi\|_A^2$ replaced by
$\int_\Omega(\varepsilon^2|\nabla\phi|^2+w\phi^2)\,\mathrm{d}x$.
\end{Lemma}
\begin{proof}
Put $a_t=(t^2-1)/2$ and $k=a_tz+th$. We first expand the quadratic part
$Q(u,v)=(u,v)_A$. Since $tz+h=(tu+\phi,tv-\phi)$,
$$
\begin{aligned}
Q(tz+h)-Q(z)
&=(t^2-1)(u,v)_A+t(\phi,v)_A-t(u,\phi)_A-\|\phi\|_A^2,\\
\langle Q'(z),k\rangle
&=(a_tu+t\phi,v)_A+(u,a_tv-t\phi)_A\\
&=(t^2-1)(u,v)_A+t(\phi,v)_A-t(u,\phi)_A.
\end{aligned}
$$
Thus
$$
Q(tz+h)-Q(z)-\langle Q'(z),k\rangle=-\|\phi\|_A^2.
$$

To treat the nonlinear part, for $s>2$ and $\xi,\psi\in\mathbb R$ define
$$
R_s(t;\xi,\psi)=-\frac{|t\xi+\psi|^s}{s}+\frac{|\xi|^s}{s}
+|\xi|^{s-2}\xi(a_t\xi+t\psi).
$$
Setting $\eta=t\xi+\psi$, we have
$$
a_t\xi+t\psi=t\eta-\frac{t^2+1}{2}\xi,
$$
and consequently
$$
R_s(t;\xi,\psi)
=t|\xi|^{s-2}\xi\eta-\frac{|\eta|^s}{s}
-\left(\frac{t^2+1}{2}-\frac1s\right)|\xi|^s.
$$
Young's inequality with conjugate exponents $s/(s-1)$ and $s$, applied to
$t|\xi|^{s-1}$ and $|\eta|$, gives
$$
\begin{aligned}
t|\xi|^{s-2}\xi\eta
&\leq t|\xi|^{s-1}|\eta|\\
&\leq\frac{s-1}{s}\bigl(t|\xi|^{s-1}\bigr)^{s/(s-1)}
+\frac{|\eta|^s}{s}\\
&=\frac{s-1}{s}t^{s/(s-1)}|\xi|^s+\frac{|\eta|^s}{s}.
\end{aligned}
$$
It follows from \eqref{Hs} that
$$
R_s(t;\xi,\psi)\leq-H_s(t)|\xi|^s.
$$
Finally, the definition of $\Phi$ and the quadratic identity above yield
$$
\begin{aligned}
&\Phi(tz+h)-\Phi(z)-\langle\Phi'(z),k\rangle\\
&\qquad=-\|\phi\|_A^2
+\int_D B R_p(t;u,\phi)\,\mathrm{d}x
+\int_D B R_q(t;v,-\phi)\,\mathrm{d}x\\
&\qquad\leq-\|\phi\|_A^2-H_p(t)\int_D B|u|^p\,\mathrm{d}x
-H_q(t)\int_D B|v|^q\,\mathrm{d}x,
\end{aligned}
$$
where the last inequality uses $B\geq0$. This proves \eqref{fiberineq}.
For $I_\varepsilon$, replace $(f,g)_A$ throughout the quadratic calculation by
$\int_\Omega(\varepsilon^2\nabla f\cdot\nabla g+wfg)\,\mathrm{d}x$ and take
$B=w$ in the nonlinear terms. The same argument proves the stated estimate.
\end{proof}

\begin{Lemma}\label{fiberminimax}
For each $z\in E_\Omega\setminus E_\Omega^-$ there is a unique point
$\widehat m(z)\in\widehat E_\Omega(z)\cap\mathcal N_\varepsilon$.
This point is the strict global maximum of $I_\varepsilon$ on the fiber, and
\begin{equation}\label{minimax}
c_\varepsilon=\inf_{z\in\mathcal N_\varepsilon}I_\varepsilon(z)
=\inf_{z\in E_\Omega\setminus E_\Omega^-}
\max_{\widehat E_\Omega(z)}I_\varepsilon.
\end{equation}
The same characterization holds for $I_{\varepsilon,D}$ on every smooth
bounded subdomain $D\subseteq\Omega$, and for \eqref{generalPhi} on every bounded
open set $D$, with the corresponding constraint and least energy level.
\end{Lemma}
\begin{proof}
We prove the assertion for $\Phi$ on a bounded open set $D$; the argument for
$I_\varepsilon$ is the same. We use the notation
$$
\begin{aligned}
\|z\|^2&=\|u\|_A^2+\|v\|_A^2,\qquad z=(u,v)\in E_D,\\
\mathcal F(z)&=\int_D B\left(\frac{|u|^p}{p}+\frac{|v|^q}{q}\right)\,\mathrm{d}x,\qquad
P_z=\int_D B|u|^p\,\mathrm{d}x,\quad Q_z=\int_D B|v|^q\,\mathrm{d}x.
\end{aligned}
$$
With respect to this equivalent Hilbert norm, the decomposition
$E_D=E_D^+\oplus E_D^-$ is orthogonal and
$$
\Phi(z)=\frac12\|z^+\|^2-\frac12\|z^-\|^2-\mathcal F(z).
$$
Let $\mathcal N$ denote the corresponding generalized Nehari set. We give a
direct proof of the characterization; see also \cite[Section 4]{Szulkin-Weth}.

\noindent\textit{Step 1: maximization on a fixed fiber.}
Fix $z\notin E_D^-$ and set $e=z^+/\|z^+\|$. Then
$$
\widehat E_D(z)=\{te+h:t\geq0,\ h\in E_D^-\},\qquad
\Phi(te+h)=\frac{t^2}{2}-\frac{\|h\|^2}{2}-\mathcal F(te+h).
$$
Since $p,q>2$, we have $\Phi(te)=t^2/2-O(t^p+t^q)>0$ for small $t>0$.
We next show that $\Phi(te+h)\to-\infty$ as $t+\|h\|\to\infty$.
Write $e=(\psi,\psi)$ and choose $\delta>0$ such that
$S=\{x\in D:|\psi(x)|\geq\delta\}$ has positive measure. For
$h=(\phi,-\phi)$ and almost every $x\in S$,
$$
\max\{|t\psi+\phi|,|t\psi-\phi|\}\geq t|\psi|\geq t\delta,
$$
which is independent of the choice of $h$.
Putting $s_0=\min\{p,q\}>2$, we obtain, for $t\delta\geq1$,
$$
\mathcal F(te+h)\geq
\frac{\kappa|S|}{\max\{p,q\}}(t\delta)^{s_0},\qquad
\Phi(te+h)\leq\frac{t^2}{2}-\frac{\|h\|^2}{2}-c_e t^{s_0}.
$$
If $t$ remains bounded and $\|h\|\to\infty$, the same conclusion follows from
$\mathcal F\geq0$. Thus every maximizing sequence on this fiber is bounded.
After passing to a subsequence, $t_n\to t$ and $h_n\rightharpoonup h$ in $E_D^-$.
Compactness of the Sobolev embeddings gives
$\mathcal F(t_ne+h_n)\to\mathcal F(te+h)$, and hence
$$
\limsup_n\Phi(t_ne+h_n)\leq\Phi(te+h).
$$
The maximum is therefore attained at a point $w=t_0e+h_0$. Its value is positive,
whereas $\Phi(h)\leq0$ on $E_D^-$, so $t_0>0$. Differentiating in the fiber
directions gives
$$
\langle\Phi'(w),e\rangle=0,\qquad
\langle\Phi'(w),h\rangle=0,\qquad h\in E_D^-,
$$
and consequently $w\in\mathcal N$.

For any $w\in\mathcal N$, Lemma~\ref{fiberestimate} yields
$$
\Phi(tw+h)-\Phi(w)
\leq-\frac12\|h\|^2-H_p(t)P_w-H_q(t)Q_w,
\qquad t\geq0,\quad h\in E_D^-.
$$
Since $P_w+Q_w>0$ and $H_s(t)>0$ for $t\neq1$, the inequality is strict unless
$t=1$ and $h=0$. Also, $w\in\widehat E_D(z)\cap\mathcal N$ implies
$\widehat E_D(w)=\widehat E_D(z)$. Thus the fiber has exactly one point in
$\mathcal N$, denoted by $\widehat m(z)$, and this point is its strict global
maximum.

\noindent\textit{Step 2: attainment of the constrained minimum.}
For $y\in E_D^+$, the Sobolev inequalities give
$$
\Phi(y)\geq\frac12\|y\|^2-C(\|y\|^p+\|y\|^q).
$$
Choose $r>0$ small enough that the right-hand side is at least some $\rho>0$
when $\|y\|=r$. Every fiber contains such a point. By Step 1,
$$
d:=\inf_{\mathcal N}\Phi\geq\rho>0.
$$
Take $z_n=(u_n,v_n)\in\mathcal N$ with $\Phi(z_n)\to d$, and abbreviate
$P_n=P_{z_n}$ and $Q_n=Q_{z_n}$. The Nehari identity gives
$$
\Phi(z_n)=\Phi(z_n)-\frac12\langle\Phi'(z_n),z_n\rangle
=\left(\frac12-\frac1p\right)P_n+\left(\frac12-\frac1q\right)Q_n,
$$
so $P_n+Q_n$ is bounded. Since $z_n^+=z_n-z_n^-$, the constraint also gives
$\langle\Phi'(z_n),z_n^+\rangle=0$. By H\"older's and Sobolev inequalities,
$$
\begin{aligned}
\|z_n^+\|^2
&=\langle\mathcal F'(z_n),z_n^+\rangle\\
&\leq C\bigl(P_n^{(p-1)/p}+Q_n^{(q-1)/q}\bigr)\|z_n^+\|
\leq C\|z_n^+\|.
\end{aligned}
$$
Moreover,
$$
\|z_n^+\|^2-\|z_n^-\|^2=P_n+Q_n\geq0.
$$
Hence $(z_n)$ is bounded in $E_D$. Passing to a subsequence, we have
$$
z_n\rightharpoonup z\text{ in }E_D,\qquad
u_n\to u\text{ in }L^p(D),\qquad v_n\to v\text{ in }L^q(D),
$$
where $z=(u,v)$. In particular,
$$
\mathcal F'(z_n)\to\mathcal F'(z)\text{ in }E_D^*,\qquad
\left(\frac12-\frac1p\right)P_z+
\left(\frac12-\frac1q\right)Q_z=d>0.
$$
The negative-direction equation passes to the limit:
$\langle\Phi'(z),h\rangle=0$ for every $h\in E_D^-$. Subtracting this equation
from that for $z_n$ and taking $h=z_n^--z^-$ yields
$$
\|z_n^--z^-\|^2
=-\langle\mathcal F'(z_n)-\mathcal F'(z),z_n^--z^-\rangle\to0.
$$
Also, $z\notin E_D^-$. Indeed, otherwise testing the negative-direction equation
with $z$ would give
$$
0=\langle\Phi'(z),z\rangle=-\|z\|^2-P_z-Q_z,
$$
contradicting the preceding positive energy identity.

By Step 1, write $\widehat m(z)=t_0z+h_0\in\mathcal N$. The positive parts of
$t_0z_n+h_0$ converge weakly, their negative parts converge strongly, and their
nonlinear integrals converge. Therefore
$$
d\leq\Phi(t_0z+h_0)
\leq\liminf_n\Phi(t_0z_n+h_0)
\leq\lim_n\Phi(z_n)=d,
$$
where the last inequality follows from the maximality of $z_n$ on its fiber.
Thus $d$ is attained.

\noindent\textit{Step 3: the minimizing point is a critical point.}
Let $w\in\mathcal N$ satisfy $\Phi(w)=d$. Since $p,q>2$, we have
$\Phi\in C^2(E_D,\mathbb R)$. Apply Lemma~\ref{fiberestimate} with
$t=1+s\alpha$ and negative direction $sh$, where $\alpha\in\mathbb R$ and
$h\in E_D^-$. For $|s|$ small, $t>0$. Taylor expansion at $s=0$, together with
$\langle\Phi'(w),\alpha w+h\rangle=0$, gives
\begin{equation}\label{eq:20261113-1621}
\begin{split}
\Phi''(w)[\alpha w+h,\alpha w+h]
&\leq-\|h\|^2-C_w\alpha^2,\\
C_w&=\frac{p-2}{p-1}P_w+\frac{q-2}{q-1}Q_w>0.
\end{split}
\end{equation}
Here we used $H_r''(1)=(r-2)/(r-1)$ for $r>2$.

Fix any $\zeta\in E_D$. We seek a curve
$$
\gamma(s)=t(s)(w+s\zeta)+h(s),\qquad
t(0)=1,\quad h(0)=0,\quad h(s)\in E_D^-,
$$
satisfying
$$
\langle\Phi'(\gamma(s)),w+s\zeta\rangle=0,\qquad
\left.\Phi'(\gamma(s))\right|_{E_D^-}=0.
$$
At $(s,t,h)=(0,1,0)$ these equations hold. Their derivative with respect to
$(t,h)$ is represented on $\mathbb R\times E_D^-$ by the bilinear form
$$
\bigl((\alpha,h),(\beta,k)\bigr)
\longmapsto\Phi''(w)[\alpha w+h,\beta w+k].
$$
The negative of this bilinear form is bounded and, by
\eqref{eq:20261113-1621}, satisfies
$$
-\Phi''(w)[\alpha w+h,\alpha w+h]\geq\|h\|^2+C_w\alpha^2.
$$
It is therefore coercive on $\mathbb R\times E_D^-$. The Lax--Milgram theorem
shows that the derivative is an isomorphism. The implicit function theorem gives
$C^1$ functions $t(s),h(s)$ for $|s|$ small. Shrinking the interval if necessary,
we have $t(s)>0$ and $(w+s\zeta)^+\neq0$, and hence $\gamma(s)\in\mathcal N$.
The minimality of $w$ now implies
$$
\begin{aligned}
0=\left.\frac{\mathrm{d}}{\mathrm{d}s}\Phi(\gamma(s))\right|_{s=0}
&=\langle\Phi'(w),\zeta+t'(0)w+h'(0)\rangle\\
&=\langle\Phi'(w),\zeta\rangle.
\end{aligned}
$$
Since $\zeta$ is arbitrary, $\Phi'(w)=0$.

Conversely, a nontrivial critical point cannot belong to $E_D^-$, since
$\langle\Phi'(z),z\rangle=-\|z\|^2-P_z-Q_z<0$ there. Every nontrivial critical
point therefore lies in $\mathcal N$. It follows that $d$ is the least energy
among all nontrivial critical points. Finally, Step 1 gives
\begin{equation}\label{eq:20261113-1708}
\inf_{z\notin E_D^-}\max_{\widehat E_D(z)}\Phi
=\inf_{z\notin E_D^-}\Phi(\widehat m(z))
=\inf_{\mathcal N}\Phi=d,
\end{equation}
because $\widehat m(z)\in\mathcal N$ and $\widehat m(w)=w$ for $w\in\mathcal N$.
For $I_{\varepsilon,D}$, take $B=w$ and use the component norm defined by
$$
\|\phi\|_{\varepsilon,D}^2
=\int_D(\varepsilon^2|\nabla\phi|^2+w\phi^2)\,\mathrm{d}x.
$$
For each fixed $\varepsilon>0$ this norm is equivalent to the $H^1_0(D)$ norm,
so all three steps apply and give \eqref{minimax}.
\end{proof}

\begin{Lemma}\label{lowerenergy}
There is $\rho>0$, independent of $\varepsilon>0$, such that
$\varepsilon^{-N}c_\varepsilon\geq\rho$.
\end{Lemma}
\begin{proof}
We use the small-sphere argument of \cite[Lemma 3.6]{Tang-fiber}, keeping the
constants independent of $\varepsilon$ after rescaling.
Set $\Omega_\varepsilon=\varepsilon^{-1}\Omega$ and
$w_\varepsilon(x)=w(\varepsilon x)$. Under the change $x=\varepsilon y$, the energy
divided by $\varepsilon^N$ becomes
\begin{equation}\label{scaledenergy}
\bar I_\varepsilon(u,v)=\int_{\Omega_\varepsilon}
(\nabla u\cdot\nabla v+w_\varepsilon uv)
-\int_{\Omega_\varepsilon}w_\varepsilon
\left(\frac{|u|^p}{p}+\frac{|v|^q}{q}\right).
\end{equation}
Its least critical value is $\varepsilon^{-N}c_\varepsilon$.
The bounds $b^{-2}\leq w_\varepsilon\leq a^{-2}$ and the Sobolev inequalities
for zero extensions give
$$
\bar I_\varepsilon(\phi,\phi)
\geq c_0\|\phi\|_{H^1}^2-C\bigl(\|\phi\|_{H^1}^p+\|\phi\|_{H^1}^q\bigr).
$$
Choose $r>0$ small enough that
$$
\rho:=c_0r^2-C(r^p+r^q)>0.
$$
Let $E_\varepsilon=H^1_0(\Omega_\varepsilon)^2$ with its standard product norm.
For any $z\in E_\varepsilon\setminus E_\varepsilon^-$, write
$z^+=(\psi,\psi)\neq0$ and set
$$
y_z=\left(\frac{r\psi}{\|\psi\|_{H^1}},
\frac{r\psi}{\|\psi\|_{H^1}}\right)
=\frac{\sqrt2\,r}{\|z^+\|_{E_\varepsilon}}z^+.
$$
Since $\|z^+\|_{E_\varepsilon}=\sqrt2\|\psi\|_{H^1}$ and
$\widehat E_{\Omega_\varepsilon}(z)=\mathbb R^+z^+\oplus E_\varepsilon^-$,
the point $y_z$ belongs to this fiber. The preceding lower bound gives
$$
\max_{\widehat E_{\Omega_\varepsilon}(z)}\bar I_\varepsilon
\geq\bar I_\varepsilon(y_z)\geq\rho.
$$
Taking the infimum over all $z\in E_\varepsilon\setminus E_\varepsilon^-$ and
using \eqref{eq:20261113-1708}, we obtain
$$
\begin{aligned}
\varepsilon^{-N}c_\varepsilon
=\inf_{z\in E_\varepsilon\setminus E_\varepsilon^-}
\max_{\widehat E_{\Omega_\varepsilon}(z)}\bar I_\varepsilon\geq\inf_{z\in E_\varepsilon\setminus E_\varepsilon^-}
\bar I_\varepsilon(y_z)\geq\rho.
\end{aligned}
$$
\end{proof}

We next prove the estimate needed for varying domains and coefficients. No
regularity or geometric bounds on these domains enter its constants.

\begin{Lemma}\label{uniformfiber}
Let $D_j\subset\mathbb R^N$ be open, and let $\Phi_j$ be \eqref{generalPhi} on
$E_j=H^1_0(D_j)^2$, with $0<\kappa\leq A_j,B_j\leq K$ uniformly in $j$.
Suppose $z_j=(u_j,v_j)\in E_j$ satisfies
$$
0<C_1\leq\Phi_j(z_j)\leq C_2,\qquad
\mu_j:=\|\Phi_j'(z_j)\|_{E_j^*}\to0.
$$
Then $\{z_j\}$ is bounded. For all sufficiently large $j$,
\begin{equation}\label{fibererror}
0\leq\sup_{\widehat E_{D_j}(z_j)}\Phi_j-\Phi_j(z_j)\leq C\mu_j^2.
\end{equation}
More precisely, for all sufficiently large $j$, all $t\geq0$ and all
$\phi\in H^1_0(D_j)$,
\begin{equation}\label{fibercoercive}
\Phi_j(tz_j+(\phi,-\phi))-\Phi_j(z_j)
\leq C\mu_j^2-c\bigl((t-1)^2+\|\phi\|_{H^1}^2\bigr).
\end{equation}
The constants $C,c>0$ in these two estimates depend only on
$N,p,q,\kappa,K,C_1,C_2$.
\end{Lemma}
\begin{proof}
\noindent\textit{Step 1: uniform boundedness.}
Write $R_j=\|u_j\|_{H^1}+\|v_j\|_{H^1}$ and set
$$
P_j=\int_{D_j}B_j|u_j|^p\,\mathrm{d}x,\qquad
Q_j=\int_{D_j}B_j|v_j|^q\,\mathrm{d}x.
$$
The identity
\begin{equation}\label{PSidentity}
\Phi_j(z_j)-\tfrac12\langle\Phi_j'(z_j),z_j\rangle
=\left(\tfrac12-\tfrac1p\right)P_j+
\left(\tfrac12-\tfrac1q\right)Q_j
\end{equation}
implies $P_j+Q_j\leq C(1+\mu_j R_j)$.
Testing the derivative with $(v_j,0)$ and $(0,u_j)$, and using the uniform
Sobolev inequalities for zero extensions, gives
$$
\begin{aligned}
\|v_j\|_{H^1}&\leq C P_j^{(p-1)/p}+C\mu_j,\\
\|u_j\|_{H^1}&\leq C Q_j^{(q-1)/q}+C\mu_j.
\end{aligned}
$$
The inequalities remain valid when a component is zero, without dividing by its
norm. Since $\sigma:=\max\{(p-1)/p, (q-1)/q\}<1$ and $\mu_j\leq1$
eventually, we deduce that
$$R_j\leq C(R_j^\sigma+1),$$
which yields $R_j\leq C$. Then by $P_j+Q_j\leq C(1+\mu_j R_j)$, we also have that $\{P_j+Q_j\}$ is bounded.  Returning to \eqref{PSidentity} and using $C_1>0$, we also obtain, for large $j$,
\begin{equation}\label{nonlinearlower}
P_j+Q_j\geq c_0>0.
\end{equation}

\noindent\textit{Step 2: the error on the whole fiber.}
By Lemma~\ref{fiberestimate}, \eqref{Hlower}, \eqref{nonlinearlower} and Step 1,
there are constants $C_*,c_1>0$ such that, for all sufficiently large $j$,
$$
\begin{aligned}
\Phi_j(tz_j+(\phi,-\phi))-\Phi_j(z_j)
\leq{}&C_*\mu_j\bigl(|t^2-1|+t\|\phi\|_{H^1}\bigr)\\
&-c_1\bigl((t-1)^2+\|\phi\|_{H^1}^2\bigr),
\end{aligned}
$$
uniformly for $t\geq0$ and $\phi\in H^1_0(D_j)$.
Set $\tau=|t-1|$ and $r=\|\phi\|_{H^1}$. Since
$$
|t^2-1|\leq\tau^2+2\tau,\qquad t\leq1+\tau,\qquad
\tau r\leq\frac12(\tau^2+r^2),
$$
we obtain
$$
\begin{aligned}
\Phi_j(tz_j+(\phi,-\phi))-\Phi_j(z_j)
&\leq C_*\mu_j(\tau^2+2\tau+r+\tau r)-c_1(\tau^2+r^2)\\
&\leq C_*\mu_j\left(\frac32\tau^2+\frac12r^2+2\tau+r\right)
-c_1(\tau^2+r^2).
\end{aligned}
$$
Choose $j$ large enough that $C_*\mu_j\leq c_1/3$. Then
$$
\Phi_j(tz_j+(\phi,-\phi))-\Phi_j(z_j)
\leq-\frac{c_1}{2}(\tau^2+r^2)+2C_*\mu_j\tau+C_*\mu_j r.
$$
Young's inequality in the form $ab\leq\delta a^2+b^2/(4\delta)$,
with $\delta=c_1/4$, gives
$$
\begin{aligned}
2C_*\mu_j\tau&\leq\frac{c_1}{4}\tau^2+\frac{4C_*^2}{c_1}\mu_j^2,\\
C_*\mu_j r&\leq\frac{c_1}{4}r^2+\frac{C_*^2}{c_1}\mu_j^2.
\end{aligned}
$$
Consequently,
$$
\Phi_j(tz_j+(\phi,-\phi))-\Phi_j(z_j)
\leq\frac{5C_*^2}{c_1}\mu_j^2
-\frac{c_1}{4}\bigl((t-1)^2+\|\phi\|_{H^1}^2\bigr).
$$
This proves \eqref{fibercoercive}, with constants independent of $t$ and $\phi$.
Since $t=1$ and $\phi=0$ are admissible, taking the supremum gives
$$
0\leq\sup_{\widehat E_{D_j}(z_j)}\Phi_j-\Phi_j(z_j)
\leq\frac{5C_*^2}{c_1}\mu_j^2,
$$
which is \eqref{fibererror}.
\end{proof}

\begin{Remark}\label{uniformremark}
The lemma includes the unweighted functional on expanding balls and the weighted
functional \eqref{scaledenergy}. It also applies after translations of the
rescaled domains. If a fiber maximum is attained, \eqref{fibercoercive} yields
$|t_j-1|+\|\phi_j\|_{H^1}=O(\mu_j)$ for its parameters. The proof is independent
of any comparison condition relating $|u|^{p-2}uv$ and $|v|^{q-2}vu$.
\end{Remark}

\section{The upper bound for the least energy}\label{sec:5}

Let $c_\infty$ be the least energy among nontrivial critical points of
$$
I_\infty(U,V)=\int_{\mathbb R^N}(\nabla U\cdot\nabla V+UV)
-\frac1p\int_{\mathbb R^N}|U|^p-\frac1q\int_{\mathbb R^N}|V|^q.
$$
The associated system is

\begin{equation}\label{5.0}
\begin{cases}
-\Delta U+U=|V|^{q-2}V &\text{in }\mathbb R^N,\\
-\Delta V+V=|U|^{p-2}U &\text{in }\mathbb R^N,\\
U,V\in H^1(\mathbb R^N).
\end{cases}
\end{equation}

\begin{Lemma}\label{autonomous}
The value $c_\infty>0$ is attained by a positive pair. Every positive
$H^1$ solution of \eqref{5.0} is, up to a common translation, radial and strictly
decreasing in the radial variable. Each such solution satisfies
$$
U(x)+V(x)+|\nabla U(x)|+|\nabla V(x)|\leq Ce^{-c|x|}.
$$
For $\lambda>0$, the least energy $c_\lambda$ of
\begin{equation}\label{frozen}
-\Delta U+\lambda U=\lambda|V|^{q-2}V,\qquad
-\Delta V+\lambda V=\lambda|U|^{p-2}U
\end{equation}
on $H^1(\mathbb R^N)^2$ is
\begin{equation}\label{frozenenergy}
c_\lambda=\lambda^{1-N/2}c_\infty.
\end{equation}
\end{Lemma}
\begin{proof}
Existence of a ground state with both components positive follows from
\cite[Theorems 1.8--1.9]{Bonheure-Ramos-ground}; see also
\cite[Theorem 6.4]{Bonheure-Tavares}. In their notation the exponents are
$p-1$ and $q-1$, and the hypotheses hold because
$1>1/p+1/q>(N-2)/N$. The generalized Nehari characterization also follows from
\cite[Theorem 42 and its proof]{Szulkin-Weth}. Its positive lower bound is obtained by the
same small-sphere estimate as in Lemma~\ref{lowerenergy}.
For a positive $H^1$ solution, local subsolution estimates imply boundedness and
decay at infinity. The nonlinearities are $C^1$, and the system is cooperative
and fully coupled on the positive quadrant. Its linearization at zero is
$-\mathrm{Id}$, so the principal minors of its negative remain positive near zero.
Thus the hypotheses of \cite[Theorem 1]{BS} hold, and both components are radial
about the same point and strictly decreasing in the radial variable.

For $S=U+V$, smallness at infinity gives
$-\Delta S+\tfrac12S\leq0$ outside a ball. Comparison with
$Ce^{-c|x|}$, where $0<c<1/\sqrt2$, gives exponential decay. Local elliptic
estimates give the same decay for the gradients. Finally,
$(U,V)\mapsto(U(\sqrt\lambda\,\cdot),V(\sqrt\lambda\,\cdot))$ is a bijection
between solutions of \eqref{5.0} and \eqref{frozen}; its energy multiplier is
$\lambda^{1-N/2}$. This proves \eqref{frozenenergy}.
\end{proof}

We also record the domain comparison, which follows directly from the best
constant formulation and does not require an ordering of indefinite quadratic
forms.
\begin{Lemma}\label{domainmonotonic}
If $D_1\subseteq D_2\subseteq\Omega$ are smooth bounded connected domains, then
$c_{\varepsilon,D_2}\leq c_{\varepsilon,D_1}$ for each fixed $\varepsilon>0$.
\end{Lemma}
\begin{proof}
Let $(u_1,v_1)$ be a positive ground state on $D_1$ and put
$K=\int_{D_1}wu_1^p\mathrm{d}x=\int_{D_1}wv_1^q\mathrm{d}x$.
Extend $wv_1^{q-1}$ by zero to $D_2$, calling the resulting function $f$, and
let $\widehat u$ solve $L_\varepsilon\widehat u=f$ in $D_2$ with zero boundary
values. The maximum principle gives $\widehat u\geq u_1$ on $D_1$ and
$\widehat u\geq0$ on $D_2$. Hence
$$
A_{D_2}(\widehat u)=K,\qquad B_{D_2}(\widehat u)\geq K,
\qquad S_{\varepsilon,D_2}\leq K^{1-q'/p}.
$$
By \eqref{energyconstant},
$c_{\varepsilon,D_2}\leq\delta_{p,q}K=c_{\varepsilon,D_1}$.
\end{proof}

\begin{Lemma}\label{upperenergy}
We have
\begin{equation}\label{upperbound}
\limsup_{\varepsilon\to0}\varepsilon^{-N}c_\varepsilon
\leq a^{N-2}c_\infty.
\end{equation}
\end{Lemma}
\begin{proof}
Fix $\xi\in\Omega$ and let $\lambda=w(\xi)$. By Lemma~\ref{autonomous},
choose a positive ground state $(U,V)$ of \eqref{frozen}. Choose
$\chi\in C_c^\infty(\Omega)$ with $0\leq\chi\leq1$ and $\chi=1$ in a
neighborhood of $\xi$. On $D_\varepsilon=(\Omega-\xi)/\varepsilon$, set
$$
U_\varepsilon(y)=\chi(\xi+\varepsilon y)U(y),\qquad
V_\varepsilon(y)=\chi(\xi+\varepsilon y)V(y),\qquad
w_\varepsilon(y)=w(\xi+\varepsilon y).
$$
After zero extension, $(U_\varepsilon,V_\varepsilon)\to(U,V)$ strongly in
$H^1(\mathbb R^N)^2$. Let $\Phi_\varepsilon$ be \eqref{generalPhi} on
$D_\varepsilon$ with $A=B=w_\varepsilon$. We claim that
\begin{equation}\label{cutoffPS}
\Phi_\varepsilon(U_\varepsilon,V_\varepsilon)=c_\lambda+o(1),
\qquad \|\Phi_\varepsilon'(U_\varepsilon,V_\varepsilon)\|\to0.
\end{equation}
To prove \eqref{cutoffPS}, extend $\chi$ by zero to $\mathbb R^N$ and write
$\chi_\varepsilon(y)=\chi(\xi+\varepsilon y)$. Choose $\delta>0$ such that
$\chi=1$ on $B_\delta(\xi)$. For $W=U,V$, set
$W_\varepsilon=\chi_\varepsilon W\in H^1_0(D_\varepsilon)$.
Since $\|\nabla\chi_\varepsilon\|_\infty\leq C\varepsilon$, we have
$$
\|W_\varepsilon-W\|_{H^1(\mathbb R^N)}
\leq C\|W\|_{H^1(\mathbb R^N\setminus B_{\delta/\varepsilon}(0))}
+C\varepsilon\|W\|_2\longrightarrow0.
$$
In particular, with $z_\varepsilon=(U_\varepsilon,V_\varepsilon)$ and $z=(U,V)$,
$$
d_\varepsilon:=\|U_\varepsilon-U\|_{H^1}
+\|V_\varepsilon-V\|_{H^1}\longrightarrow0.
$$

Put $g_\varepsilon=w_\varepsilon-\lambda$ on $D_\varepsilon$. For each fixed
$R>0$, we have $B_R(0)\subset D_\varepsilon$ for small $\varepsilon$, and
$$
\omega_\varepsilon(R):=
\|g_\varepsilon\|_{L^\infty(B_R(0))}\longrightarrow0,
\qquad \|g_\varepsilon\|_{L^\infty(D_\varepsilon)}\leq C.
$$
For $(W,s)=(U,2),(V,2),(U,p),(V,q)$, let $s'=s/(s-1)$ and set
$$
\mathcal E_\varepsilon(W,s)
:=\bigl\|g_\varepsilon|W_\varepsilon|^{s-1}
\bigr\|_{L^{s'}(D_\varepsilon)}.
$$
Using $|W_\varepsilon|\leq|W|$ and splitting the norm over $B_R(0)$ and its
complement, we obtain
$$
\mathcal E_\varepsilon(W,s)
\leq\omega_\varepsilon(R)\|W\|_s^{s-1}
+C\|W\|_{L^s(\mathbb R^N\setminus B_R(0))}^{s-1}.
$$
First letting $\varepsilon\to0$ and then $R\to\infty$ gives
$\mathcal E_\varepsilon(W,s)\to0$ for all four pairs.

Let $I_\lambda$ denote the energy of \eqref{frozen} on $H^1(\mathbb R^N)^2$. By
H\"older's inequality and the Sobolev embeddings,
$|I_\lambda(z_\varepsilon)-I_\lambda(z)|\leq C d_\varepsilon$.
Since $I_\lambda(z)=c_\lambda$ and
$$
\Phi_\varepsilon(z_\varepsilon)-I_\lambda(z_\varepsilon)
=\int_{D_\varepsilon}g_\varepsilon
\left(U_\varepsilon V_\varepsilon
-\frac{|U_\varepsilon|^p}{p}-\frac{|V_\varepsilon|^q}{q}\right)\,dy,
$$
we conclude that
$$
\begin{aligned}
|\Phi_\varepsilon(z_\varepsilon)-c_\lambda|
\leq{}&C d_\varepsilon
+\mathcal E_\varepsilon(U,2)\|V_\varepsilon\|_2\\
&+\frac1p\mathcal E_\varepsilon(U,p)\|U_\varepsilon\|_p
+\frac1q\mathcal E_\varepsilon(V,q)\|V_\varepsilon\|_q
\longrightarrow0.
\end{aligned}
$$

For the derivative, put $F_s(a)=|a|^{s-2}a$. If $(W,s)=(U,p)$ or $(V,q)$,
then
$$
\begin{aligned}
\|F_s(W_\varepsilon)-F_s(W)\|_{s'}
&\leq C_s(\|W_\varepsilon\|_s+\|W\|_s)^{s-2}
\|W_\varepsilon-W\|_s\\
&\leq C d_\varepsilon\longrightarrow0.
\end{aligned}
$$
Take $(\phi,\psi)\in H^1_0(D_\varepsilon)^2$ and extend both components by
zero. Since $I_\lambda'(z)=0$, we obtain
$$
\begin{aligned}
\langle\Phi_\varepsilon'(z_\varepsilon),(\phi,\psi)\rangle
={}&\int_{\mathbb R^N}
\bigl[\nabla(V_\varepsilon-V)\cdot\nabla\phi
+\nabla(U_\varepsilon-U)\cdot\nabla\psi\bigr]\\
&+\lambda\int_{\mathbb R^N}
\bigl[(V_\varepsilon-V)\phi+(U_\varepsilon-U)\psi\bigr]\\
&-\lambda\int_{\mathbb R^N}
\bigl[(F_p(U_\varepsilon)-F_p(U))\phi
+(F_q(V_\varepsilon)-F_q(V))\psi\bigr]\\
&+\int_{D_\varepsilon}g_\varepsilon
\bigl[V_\varepsilon\phi+U_\varepsilon\psi
-F_p(U_\varepsilon)\phi-F_q(V_\varepsilon)\psi\bigr].
\end{aligned}
$$
The Sobolev inequalities for the zero extensions give
$$
\|\phi\|_p+\|\psi\|_q
\leq C\|(\phi,\psi)\|_{H^1_0(D_\varepsilon)^2},
$$
with $C$ independent of $\varepsilon$. Hence the preceding estimates yield
$$
\begin{aligned}
\|\Phi_\varepsilon'(z_\varepsilon)\|_{(H^1_0(D_\varepsilon)^2)^*}
\leq C\bigl(&d_\varepsilon
+\mathcal E_\varepsilon(U,2)+\mathcal E_\varepsilon(V,2)\\
&+\mathcal E_\varepsilon(U,p)+\mathcal E_\varepsilon(V,q)\bigr)
\longrightarrow0.
\end{aligned}
$$
This proves \eqref{cutoffPS}.

The least critical value of $\Phi_\varepsilon$ is
$\varepsilon^{-N}c_\varepsilon$. Since $c_\lambda>0$,
Lemma~\ref{uniformfiber} and the minimax characterization give
$$
\varepsilon^{-N}c_\varepsilon
\leq\max_{\widehat E_{D_\varepsilon}(U_\varepsilon,V_\varepsilon)}
\Phi_\varepsilon
\leq\Phi_\varepsilon(U_\varepsilon,V_\varepsilon)+o(1)
=c_\lambda+o(1).
$$
By \eqref{frozenenergy}, $c_\lambda=(2|\xi|)^{(N-2)/2}c_\infty$.
First let $\varepsilon\to0$ with $\xi$ fixed, and then let $|\xi|\downarrow a^2/2$.
This proves \eqref{upperbound}.
\end{proof}

\section{Concentration of solutions}\label{sec:6}

For each $\varepsilon>0$, choose a positive ground state $(u_\varepsilon,v_\varepsilon)$
of \eqref{2.3}, provided by Lemma~\ref{groundinversion}. By
Proposition~\ref{commonsymmetry}, the two components have a common foliated
Schwarz symmetry direction. We first analyze their concentration without fixing
this direction.

\subsection{A half-space nonexistence result}
\begin{Lemma}\label{halfspace}
Let $p,q>2$ and $\lambda,\mu_1,\mu_2>0$. Let $\mathbb R^N_+=\{y\in\mathbb R^N:y_N>0\}$. Suppose that
$u,v\in C^2(\mathbb R^N_+)\cap C^1(\overline{\mathbb R^N_+})$ are nonnegative,
vanish at infinity, and satisfy
\begin{equation}\label{6.1}
\begin{cases}
-\Delta u+\lambda u=\mu_1v^{q-1} &\text{in }\mathbb R^N_+,\\
-\Delta v+\lambda v=\mu_2u^{p-1} &\text{in }\mathbb R^N_+,\\
u=v=0 &\text{on }\partial\mathbb R^N_+.
\end{cases}
\end{equation}
Then $u=v=0$.
\end{Lemma}
\begin{proof}
Suppose that the pair is nontrivial. The strong maximum principle gives $u,v>0$
in the half-space, and Hopf's lemma gives
\begin{equation}\label{hopf}
u_{y_N}(y',0)>0,\qquad v_{y_N}(y',0)>0.
\end{equation}
We first justify integration over the unbounded domain. Since $u,v\to0$ at
infinity, $S=u+v$ satisfies
$-\Delta S+(\lambda/2)S\leq0$ outside a sufficiently large ball. Comparison
with $Ce^{-c|y|}$, $0<c<\sqrt{\lambda/2}$, gives $S\leq Ce^{-c|y|}$.
The comparison is valid on the flat boundary because $S=0$, and at infinity
because both functions tend to zero. Interior and flat-boundary elliptic
estimates then give the same exponential bound for $|\nabla u|+|\nabla v|$.

Multiply the first equation by $v_{y_N}$ and the second by $u_{y_N}$, add them,
and integrate on $B_R(0)\cap\mathbb R^N_+$. The terms without Laplacians are
the $y_N$ derivative of
$\lambda uv-\mu_1v^q/q-\mu_2u^p/p$.
The spherical boundary terms tend to zero as $R\to\infty$ by the exponential
estimates. On the flat boundary the functions and their tangential derivatives
vanish. Integration by parts consequently yields
$$
0=\int_{\mathbb R^N_+}(\Delta u\,v_{y_N}+\Delta v\,u_{y_N})\,dy
=-\int_{\mathbb R^{N-1}}u_{y_N}(y',0)v_{y_N}(y',0)\,dy'.
$$
The last integral is finite and strictly positive by \eqref{hopf}, a
contradiction. This is the translation identity used for scalar and weakly
coupled equations in \cite{Estaban-Lions,Lin-Wei}.
\end{proof}

\subsection{Uniform bounds and exclusion of boundary profiles}
Let $Z_\varepsilon$ be any maximum point of $u_\varepsilon+v_\varepsilon$.
Such a point belongs to $\Omega$ by positivity and the Dirichlet condition.

\begin{Lemma}\label{uniformregularity}
For any $\xi_\varepsilon\in\overline\Omega$, define
$$
D_\varepsilon=(\Omega-\xi_\varepsilon)/\varepsilon,\qquad
U_\varepsilon(y)=u_\varepsilon(\xi_\varepsilon+\varepsilon y),\quad
V_\varepsilon(y)=v_\varepsilon(\xi_\varepsilon+\varepsilon y).
$$
Their zero extensions satisfy
\begin{equation}\label{uniformHV}
\|U_\varepsilon\|_{H^1(\mathbb R^N)}+
\|V_\varepsilon\|_{H^1(\mathbb R^N)}+
\|U_\varepsilon\|_\infty+\|V_\varepsilon\|_\infty\leq C.
\end{equation}
If $W_\varepsilon$ is the zero extension of $U_\varepsilon+V_\varepsilon$, then
for every $y\in\mathbb R^N$,
\begin{equation}\label{globalball}
\sup_{B_1(y)}W_\varepsilon
\leq C\left(\int_{B_2(y)}W_\varepsilon^2\right)^{1/2}.
\end{equation}
Furthermore, $\|u_\varepsilon\|_\infty\geq1$ and
$\|v_\varepsilon\|_\infty\geq1$.
\end{Lemma}
\begin{proof}
By Lemmas~\ref{lowerenergy} and \ref{upperenergy}, the rescaled energies lie in
a fixed positive bounded interval. The rescaled pair is a critical point of
\eqref{generalPhi} with $A=B=w(\xi_\varepsilon+\varepsilon y)$.
Step 1 of Lemma~\ref{uniformfiber} gives the $H^1$ bounds.

We explain why the local bound also applies to balls meeting the boundary.
For each fixed $\varepsilon$, the pair is classical and positive with zero
Dirichlet values. Its outward normal derivatives are nonpositive. Green's
formula therefore gives, in the sense of distributions on all of $\mathbb R^N$,
\begin{equation}\label{zerosubsolution}
-\Delta W_\varepsilon
\leq C\bigl(W_\varepsilon^{p-1}+W_\varepsilon^{q-1}\bigr)
\leq C(1+W_\varepsilon^{s-2})W_\varepsilon,
\qquad s=\max\{p,q\}<2^*.
\end{equation}
In fact the boundary contribution is
$\partial_\nu(U_\varepsilon+V_\varepsilon)$ times the surface measure, and has
the required nonpositive sign.
Choose $t>N/2$ with $(s-2)t\leq2^*$. The $H^1$ bound shows that
$C(1+W_\varepsilon^{s-2})$ has uniformly bounded $L^t$ norm on every ball
of radius two. The local estimate for nonnegative subsolutions
\cite[Theorem 8.17]{GT} applied to \eqref{zerosubsolution} proves
\eqref{globalball}, with a constant independent of the center and of
$\varepsilon$. This also proves the $L^\infty$ bounds in \eqref{uniformHV}.

Finally, let $A_\varepsilon=\max u_\varepsilon$ and
$B_\varepsilon=\max v_\varepsilon$. Evaluating the two equations at their
respective maximum points gives
\begin{equation}\label{maxineq}
A_\varepsilon\leq B_\varepsilon^{q-1},\qquad
B_\varepsilon\leq A_\varepsilon^{p-1}.
\end{equation}
Since $A_\varepsilon,B_\varepsilon>0$ and $(p-1)(q-1)>1$, these inequalities
imply $A_\varepsilon,B_\varepsilon\geq1$.
\end{proof}

\begin{Lemma}\label{boundarydistance}
As $\varepsilon\to0$,
\begin{equation}\label{Zdistance}
\frac{\operatorname{dist}(Z_\varepsilon,\partial\Omega)}{\varepsilon}\to\infty.
\end{equation}
\end{Lemma}
\begin{proof}
Suppose by contradiction that there are $\varepsilon_j\to0$ and $M>0$ such that,
with $Z_j=Z_{\varepsilon_j}$,
$$
d_j:=\operatorname{dist}(Z_j,\partial\Omega),\qquad
t_j:=\frac{d_j}{\varepsilon_j}\leq M.
$$
We shall obtain a nontrivial half-space solution, contrary to Lemma~\ref{halfspace}.

\noindent\textit{Step 1: the half-space limit.}
Choose a nearest boundary point $\zeta_j$. Since $d_j\leq M\varepsilon_j\to0$,
the tubular neighborhood property of the fixed smooth boundary gives, for large $j$,
$$
Z_j=\zeta_j+d_j\nu_{\rm in}(\zeta_j),
$$
where $\nu_{\rm in}$ denotes the inward unit normal. Let $R_j$ be orthogonal with
$R_je_N=\nu_{\rm in}(\zeta_j)$, and set
$$
\begin{aligned}
D_j&=R_j^{-1}(\Omega-\zeta_j)/\varepsilon_j,\qquad
w_j(y)=w(\zeta_j+\varepsilon_jR_jy),\\
U_j(y)&=u_{\varepsilon_j}(\zeta_j+\varepsilon_jR_jy),\qquad
V_j(y)=v_{\varepsilon_j}(\zeta_j+\varepsilon_jR_jy).
\end{aligned}
$$
Thus the maximum point becomes $y_j=t_je_N$, and
$$
\begin{cases}
-\Delta U_j+w_jU_j=w_jV_j^{q-1}&\text{in }D_j,\\
-\Delta V_j+w_jV_j=w_jU_j^{p-1}&\text{in }D_j,\\
U_j=V_j=0&\text{on }\partial D_j.
\end{cases}
$$

In coordinates centered at $\zeta_j$ with inward normal $e_N$, the original
boundary has a graph $x_N=h_j(x')$ on a fixed neighborhood, with
$h_j(0)=0$, $\nabla h_j(0)=0$, and uniformly bounded $C^2$ norm.
The domain $D_j$ therefore lies locally above the graph
$$
g_j(y')=\frac{h_j(\varepsilon_jy')}{\varepsilon_j},\qquad
|g_j(y')|\leq C\varepsilon_j|y'|^2,\qquad
|\nabla g_j(y')|\leq C\varepsilon_j|y'|.
$$
Also, $\|D^2g_j\|_\infty\leq C\varepsilon_j$. Hence $D_j$ converges locally
to $\mathbb R^N_+$, with uniform $C^{1,1}$ boundary charts on each fixed ball.
Passing to a subsequence, $\zeta_j\to\zeta\in\partial\Omega$ and
$w_j\to\lambda:=w(\zeta)>0$ uniformly on compact subsets of $\mathbb R^N_+$.

Lemma~\ref{uniformregularity}, which is unchanged by rotations, gives uniform
$H^1$ and $L^\infty$ bounds for the zero extensions of $U_j,V_j$.
Weak compactness and interior elliptic estimates yield, along a subsequence,
$$
\begin{aligned}
(U_j,V_j)&\rightharpoonup(U,V)\quad\text{in }H^1(\mathbb R^N)^2,\\
(U_j,V_j)&\longrightarrow(U,V)\quad\text{in }
C^1_{\rm loc}(\mathbb R^N_+)^2.
\end{aligned}
$$
Here the second convergence concerns the unextended functions on compact subsets
eventually contained in $D_j$. Every compact subset of $\mathbb R^N_-$ is
eventually outside $D_j$, so $U=V=0$ there. Since the limits belong to
$H^1(\mathbb R^N)$, their traces from the two sides of $y_N=0$ agree.
Consequently, their restrictions belong to $H^1_0(\mathbb R^N_+)$.
Passing to the limit in the interior equations shows that $U,V$ are bounded,
nonnegative solutions of \eqref{6.1} with $\mu_1=\mu_2=\lambda$.

\noindent\textit{Step 2: the limit is nontrivial.}
By Lemma~\ref{uniformregularity} and the choice of $Z_j$,
$$
(U_j+V_j)(t_je_N)
=\max_\Omega(u_{\varepsilon_j}+v_{\varepsilon_j})\geq1.
$$
After taking a subsequence, $t_j\to t_0\in[0,M]$. We first rule out $t_0=0$.
The uniform $L^\infty$ bounds imply that
$w_j(V_j^{q-1}-U_j)$ and $w_j(U_j^{p-1}-V_j)$ are uniformly bounded.
Local Dirichlet estimates on $D_j\cap B_2(0)$, with zero data on the true
boundary $\partial D_j$, give, for a fixed $\ell>N$,
$$
\|U_j\|_{W^{2,\ell}(D_j\cap B_1(0))}
+\|V_j\|_{W^{2,\ell}(D_j\cap B_1(0))}\leq C.
$$
The uniform boundary charts and Sobolev embedding thus give
$|\nabla U_j|+|\nabla V_j|\leq C$ near the origin.
If $t_j\to0$, the segment $\{se_N:0<s\leq t_j\}$ lies in this neighborhood
inside $D_j$. Using $U_j(0)=V_j(0)=0$, we obtain
$$
1\leq(U_j+V_j)(t_je_N)
=\int_0^{t_j}\partial_{y_N}(U_j+V_j)(se_N)\,ds
\leq Ct_j\longrightarrow0,
$$
a contradiction. Therefore $t_0>0$. The points $t_je_N$ now lie in a fixed
compact subset of $\mathbb R^N_+$, and the local uniform convergence gives
$$
(U+V)(t_0e_N)=\lim_{j\to\infty}(U_j+V_j)(t_je_N)\geq1.
$$

\noindent\textit{Step 3: decay and contradiction.}
Interior and flat-boundary regularity give
$U,V\in C^2(\mathbb R^N_+)\cap C^1(\overline{\mathbb R^N_+})$, with zero
boundary values. Let $W$ be the zero extension of $U+V$.
Boundedness of $U,V$ and the nonpositive outward boundary derivative imply,
as in the proof of Lemma~\ref{uniformregularity},
$-\Delta W\leq CW$ in the sense of distributions on $\mathbb R^N$.
The local subsolution estimate therefore gives, for $|y|>2$,
$$
\begin{aligned}
0\leq W(y)&\leq\sup_{B_1(y)}W
\leq C\left(\int_{B_2(y)}W^2\,dx\right)^{1/2}\\
&\leq C\left(\int_{\{|x|>|y|-2\}}W^2\,dx\right)^{1/2}
\longrightarrow0\qquad\text{as }|y|\to\infty,
\end{aligned}
$$
where $C$ is independent of $y$ and the last limit follows from
$W\in L^2(\mathbb R^N)$. Thus $U,V$ vanish at infinity.
Lemma~\ref{halfspace} now implies $U=V=0$, contradicting Step 2.
This proves \eqref{Zdistance}.
\end{proof}

\subsection{The limiting energy and the concentration profile}
\begin{Lemma}\label{profile}
We have
\begin{equation}\label{energylimit}
|Z_\varepsilon|\to\frac{a^2}{2},\qquad
\varepsilon^{-N}c_\varepsilon\to a^{N-2}c_\infty.
\end{equation}
There are $C,c>0$ such that
\begin{equation}\label{sumdecay}
u_\varepsilon(x)+v_\varepsilon(x)
\leq Ce^{-c|x-Z_\varepsilon|/\varepsilon}\qquad(x\in\Omega).
\end{equation}
Moreover, for every sequence $\varepsilon_j\to0$, a subsequence of the profiles
$$
U_j(y)=u_{\varepsilon_j}(Z_{\varepsilon_j}+\varepsilon_j y),\qquad
V_j(y)=v_{\varepsilon_j}(Z_{\varepsilon_j}+\varepsilon_j y)
$$
converges, after zero extension, strongly in $H^1(\mathbb R^N)^2$ to a positive
ground state of \eqref{frozen} with $\lambda=a^{-2}$. On every fixed compact
subset of $\mathbb R^N$, the unextended profiles are defined for all sufficiently
large $j$ and converge in $C^2$. Both limiting components are radial about zero
and strictly decreasing in $|y|$.
\end{Lemma}
\begin{proof}
\noindent\textit{Step 1: equality of the limiting energy bounds.}
Take an arbitrary sequence $\varepsilon_j\to0$ and abbreviate
$Z_j=Z_{\varepsilon_j}$. By compactness of $\overline\Omega$, we may assume
$Z_j\to Z_0$ and set $d=|Z_0|\geq a^2/2$, $\lambda=(2d)^{-1}$.
By Lemma~\ref{boundarydistance}, the domains
$D_j=(\Omega-Z_j)/\varepsilon_j$ exhaust $\mathbb R^N$ locally.
Lemma~\ref{uniformregularity} and elliptic estimates yield, along a subsequence,
$$
(U_j,V_j)\rightharpoonup(U,V)\text{ in }H^1(\mathbb R^N)^2,
\qquad (U_j,V_j)\to(U,V)\text{ in }C^2_{\rm loc}(\mathbb R^N)^2.
$$
The limit solves \eqref{frozen}, is nonnegative, and satisfies
$(U+V)(0)\geq1$. Both components are strictly positive by the maximum principle.

Let $\alpha_p=1/2-1/p$ and $\alpha_q=1/2-1/q$, which are positive.
The critical point identity gives
\begin{equation}\label{energydensity}
\varepsilon_j^{-N}c_{\varepsilon_j}
=\int_{D_j}w(Z_j+\varepsilon_jy)
\bigl(\alpha_p U_j^p+\alpha_q V_j^q\bigr)\,dy.
\end{equation}
Fatou's lemma, \eqref{frozenenergy} and Lemma~\ref{upperenergy} imply
\begin{equation}\label{sandwich}
\begin{split}
a^{N-2}c_\infty
&\geq\limsup_j\varepsilon_j^{-N}c_{\varepsilon_j}
\geq\liminf_j\varepsilon_j^{-N}c_{\varepsilon_j}\\
&\geq\lambda\int_{\mathbb R^N}(\alpha_p U^p+\alpha_q V^q)
=I_\lambda(U,V)\geq c_\lambda
=(2d)^{(N-2)/2}c_\infty.
\end{split}
\end{equation}
Here $I_\lambda$ is the energy associated with \eqref{frozen}.
Since $N=m+1\geq3$ and $d\geq a^2/2$, all inequalities are equalities.
Thus $d=a^2/2$ and $(U,V)$ is a ground state for $\lambda=a^{-2}$.
Since every initial sequence admits this conclusion, \eqref{energylimit} holds
for the full family.

\noindent\textit{Step 2: strong convergence.}
Extend $w_j(y)=w(Z_j+\varepsilon_jy)$ from $D_j$ to all of $\mathbb R^N$ by
the value $\lambda$. Then $b^{-2}\leq w_j\leq a^{-2}$ and $w_j\to\lambda$
pointwise. Equality in \eqref{sandwich}, together with Fatou's lemma for each
of its two nonnegative terms, gives
\begin{equation}\label{weightednorms}
\int_{\mathbb R^N}w_jU_j^p\to\lambda\int_{\mathbb R^N}U^p,
\qquad
\int_{\mathbb R^N}w_jV_j^q\to\lambda\int_{\mathbb R^N}V^q.
\end{equation}
The Br\'ezis--Lieb lemma \cite{Brezis-Lieb} applied to $w_j^{1/p}U_j$ and $w_j^{1/q}V_j$,
and dominated convergence applied to the fixed limits, yield
\begin{equation}\label{strongLp}
U_j\to U\text{ in }L^p(\mathbb R^N),\qquad
V_j\to V\text{ in }L^q(\mathbb R^N).
\end{equation}
In particular,
$w_jV_j^{q-1}\to\lambda V^{q-1}$ in $L^{q'}$ and
$w_jU_j^{p-1}\to\lambda U^{p-1}$ in $L^{p'}$.
Testing the first equation with $U_j$, and using the weak convergence of $U_j$
in $L^q$ as well as the first of these strong convergences, gives
\begin{equation}\label{quadraticlimit}
\int_{\mathbb R^N}(|\nabla U_j|^2+w_jU_j^2)
=\int_{\mathbb R^N}w_jV_j^{q-1}U_j
\longrightarrow\int_{\mathbb R^N}(|\nabla U|^2+\lambda U^2).
\end{equation}
Also, $w_jU\to\lambda U$ in $L^2$. Expanding the quadratic expression for
$U_j-U$, using weak convergence and \eqref{quadraticlimit}, we obtain
$$
\int_{\mathbb R^N}\bigl(|\nabla(U_j-U)|^2+w_j|U_j-U|^2\bigr)\to0.
$$
The uniform positive lower bound on $w_j$ proves strong $H^1$ convergence.
Here weak $L^q$ convergence follows from weak $H^1$ convergence and the
continuous embedding $H^1(\mathbb R^N)\hookrightarrow L^q(\mathbb R^N)$.
For $V_j$, the same argument uses weak $L^p$ convergence.

By Lemma~\ref{autonomous} and scaling, $U,V$ are strictly radially decreasing
about a common point. Since $U_j+V_j$ has a global maximum at zero, so does
$U+V$. The common center is therefore zero.

\noindent\textit{Step 3: uniform decay, including near the boundary.}
Let $W_\varepsilon$ be the zero extension of the sum of the profiles centered
at $Z_\varepsilon$. We claim that
\begin{equation}\label{uniformtail}
\lim_{R\to\infty}\limsup_{\varepsilon\to0}
\sup_{|y|\geq R}W_\varepsilon(y)=0.
\end{equation}
If not, there would be $\varepsilon_j\to0$ and $|y_j|\to\infty$ with
$W_{\varepsilon_j}(y_j)\geq\eta>0$. Step 2 gives, after a subsequence, strong
$L^2$ convergence of these centered profiles. Hence
$\int_{B_2(y_j)}W_{\varepsilon_j}^2\to0$. This contradicts
\eqref{globalball}. No containment of $B_2(y_j)$ in the rescaled domain is
needed here.

By \eqref{uniformtail}, choose $R$ and $\varepsilon_0>0$ such that, on
$D_\varepsilon\setminus B_R(0)$ and for $0<\varepsilon<\varepsilon_0$,
$$
U_\varepsilon^{p-1}+V_\varepsilon^{q-1}
\leq\tfrac12(U_\varepsilon+V_\varepsilon).
$$
Thus
$$
-\Delta(U_\varepsilon+V_\varepsilon)
+\frac1{2b^2}(U_\varepsilon+V_\varepsilon)\leq0
\quad\text{on }D_\varepsilon\setminus B_R(0).
$$
Choose $0<c<1/(\sqrt2 b)$ and compare with $M e^{-c(|y|-R)}$, where $M$
is a uniform bound from Lemma~\ref{uniformregularity}. This is a supersolution
outside $B_R$; it dominates the sum on $|y|=R$ and on $\partial D_\varepsilon$.
The maximum principle gives $U_\varepsilon+V_\varepsilon\leq Ce^{-c|y|}$ and then
\eqref{sumdecay} follows.
\end{proof}

\subsection{The two components and the lifted solution}
An $o(\varepsilon)$ separation of the component peaks also appears in the
compact manifold problem of Kannoth and Manna \cite{Kannoth-Manna}. For the
present problem, it follows from Lemma~\ref{profile} and the strict radial
monotonicity of the limiting pair.

Choose any maximum points $P_\varepsilon$ and $Q_\varepsilon$ of
$u_\varepsilon$ and $v_\varepsilon$, respectively.
\begin{Lemma}\label{componentpeaks}
As $\varepsilon\to0$,
\begin{gather}
|P_\varepsilon|,|Q_\varepsilon|\to a^2/2,
\label{peakradii}\\
\frac{\operatorname{dist}(P_\varepsilon,\partial\Omega)}{\varepsilon},\quad
\frac{\operatorname{dist}(Q_\varepsilon,\partial\Omega)}{\varepsilon}\to\infty,
\label{peakdistance}\\
\frac{|P_\varepsilon-Z_\varepsilon|+|Q_\varepsilon-Z_\varepsilon|}
{\varepsilon}\to0.
\label{peakcoalescence}
\end{gather}
Moreover,
\begin{equation}\label{componentdecay}
u_\varepsilon(x)\leq Ce^{-c|x-P_\varepsilon|/\varepsilon},\qquad
v_\varepsilon(x)\leq Ce^{-c|x-Q_\varepsilon|/\varepsilon},
\qquad x\in\Omega.
\end{equation}
\end{Lemma}
\begin{proof}
By Lemma~\ref{uniformregularity} and \eqref{sumdecay},
$$
1\leq u_\varepsilon(P_\varepsilon)
\leq Ce^{-c|P_\varepsilon-Z_\varepsilon|/\varepsilon},\qquad
1\leq v_\varepsilon(Q_\varepsilon)
\leq Ce^{-c|Q_\varepsilon-Z_\varepsilon|/\varepsilon}.
$$
Both scaled distances are bounded. Equations \eqref{energylimit} and
\eqref{Zdistance} immediately give \eqref{peakradii} and \eqref{peakdistance}.

Take any sequence $\varepsilon_j\to0$. Along a subsequence, the profiles centered
at $Z_{\varepsilon_j}$ converge as in Lemma~\ref{profile}, and
$$
\frac{P_{\varepsilon_j}-Z_{\varepsilon_j}}{\varepsilon_j}\to P_0,\qquad
\frac{Q_{\varepsilon_j}-Z_{\varepsilon_j}}{\varepsilon_j}\to Q_0.
$$
Local uniform convergence and maximality imply that $P_0$ maximizes $U$ and
$Q_0$ maximizes $V$. The limits are strictly radially decreasing about zero,
so $P_0=Q_0=0$. The arbitrariness of the sequence proves
\eqref{peakcoalescence}; in particular $|P_\varepsilon-Q_\varepsilon|=o(\varepsilon)$.
Finally, the triangle inequality and \eqref{sumdecay} give
\eqref{componentdecay}, with uniform constants.
\end{proof}

\begin{proof}[Proof of Theorem~\ref{th1.1}]
For each $\varepsilon$, let $e_\varepsilon$ be the common direction supplied by
Proposition~\ref{commonsymmetry}. Because the annulus and its coefficients are
rotation invariant, rotate both components simultaneously so that
$e_\varepsilon=e_N$. This is a choice of the solution family for each
$\varepsilon$, and does not require convergence of the original directions.
Foliated Schwarz symmetry allows us to choose
$$
P_\varepsilon=d_\varepsilon^u e_N,\qquad
Q_\varepsilon=d_\varepsilon^v e_N
$$
as maximum points: at every fixed radius the angular maximum is attained at
polar angle zero. 

Lift the pair by \eqref{2.7} and denote the resulting solution on $A$ by
$(\widetilde u_\varepsilon,\widetilde v_\varepsilon)$. Set
$r_\varepsilon^u=\sqrt{2d_\varepsilon^u}$ and
$r_\varepsilon^v=\sqrt{2d_\varepsilon^v}$. Equation
\eqref{spherecorrespondence} gives the two spheres of maxima stated in the
theorem. The limits in \eqref{radii} follow from
\eqref{peakradii}--\eqref{peakcoalescence}; for example,
$$
\frac{r_\varepsilon^u-a}{\varepsilon}
=\frac{2}{r_\varepsilon^u+a}\,
\frac{d_\varepsilon^u-a^2/2}{\varepsilon}\to\infty.
$$
The lower bounds on the maxima are unchanged by the lift.

We verify the quantitative decay. For $x=(y_1,y_2)\in A$, let $z$ be any
point of reduced polar coordinates $(\rho,\varphi)=(|x|^2/2,2\theta)$.
For $r_0\in(a,b)$ and $d_0=r_0^2/2$, a direct calculation gives
\begin{equation}\label{distanceidentity}
\begin{split}
|z-d_0e_N|^2
&=\frac14\bigl(|x|^2+r_0^2-2r_0|y_1|\bigr)
\bigl(|x|^2+r_0^2+2r_0|y_1|\bigr)\\
&\geq\frac{a^2}{2}\,
\operatorname{dist}\bigl(x,\{(\eta,0):|\eta|=r_0\}\bigr)^2.
\end{split}
\end{equation}
Apply \eqref{componentdecay} and \eqref{distanceidentity} with
$r_0=r_\varepsilon^u$ and $r_0=r_\varepsilon^v$ to obtain \eqref{liftdecay}.
Since both radii tend to $a$, the two spheres converge to
$\Gamma=\{(y_1,0):|y_1|=a\}$ and the asserted uniform convergence away from
$\Gamma$ follows. Relabeling the lifted functions as $u_\varepsilon,v_\varepsilon$
completes the proof.
\end{proof}

\noindent\textbf{Acknowledgments.} Hui Zhang was supported by the National Social Science Fund of China (No.~25BJL036) and the Talent Program of Jinling
Institute of Technology. Jianjun Zhang was supported by the National Natural Science Foundation of China (Nos.~12371109, 11871123) and the Natural Science Foundation of Fujian Province (No.~2025J01351). Xuexiu Zhong was supported by the National Natural Science Foundation of China (No.~12271184), the Young Top-notch Talent Project of Guangdong Province (2024TQ08A725), and the Guangzhou Basic and Applied Basic Research Foundation (2024A04J10001). The authors
acknowledge the use of AI tools. All mathematical arguments and proofs in the final
manuscript were checked and written by the authors.

\end{document}